\documentclass{amsart}
\usepackage{amsmath,amssymb}
\usepackage{paralist}
\usepackage[misc]{ifsym}
\usepackage{epsfig} %For pictures: screened artwork should be set up with an 85 or 100 line screen
\usepackage{epstopdf} %This is to transfer .eps figure to .pdf figure; please compile your article using PDFLaTex or PDFTeXify.
\usepackage[colorlinks=true]{hyperref}
\hypersetup{urlcolor=blue, citecolor=red}
\allowdisplaybreaks

\newtheorem{theorem}{Theorem}[section]

\newtheorem{proposition}[theorem]{Proposition}

\theoremstyle{definition}

\newcommand{\sref}[1]{Section~\ref{#1}}

\newcommand{\Sref}[1]{Section~\ref{#1}}

\newcommand{\fref}[1]{Figure~\ref{#1}}

\newcommand{\Fref}[1]{Figure~\ref{#1}}

\newcommand{\T}{\mathbb{T}}
\newcommand{\R}{\mathbb{R}}
\newcommand{\Sp}{\mathrm{S}}
\newcommand{\Dp}{\mathrm{D}}
\newcommand{\DSS}{\mathrm{DSS}}
\newcommand{\DDS}{\mathrm{DDS}}
\newcommand{\DSD}{\mathrm{DSD}}
\newcommand{\DDD}{\mathrm{DDD}}
\newcommand{\SSS}{\mathrm{SSS}}
\newcommand{\SDS}{\mathrm{SDS}}
\newcommand{\SSD}{\mathrm{SSD}}
\newcommand{\SDD}{\mathrm{SDD}}

\newcommand{\lb}{{\boldsymbol{\ell}}}

\newcommand{\trp}{\mathsf{T}}

\title[Dynamics on invariant tori through forced symmetry breaking]
{Dynamics on invariant tori emerging through forced symmetry breaking in phase oscillator networks} % Only the first word and proper nouns should be capitalized.

\author[Christian Bick, Jos\'e Mujica, Bob Rink]{}


\thanks{JM and CB acknowledge support from the Engineering and Physical Sciences Research Council (EPSRC) through the grant EP/T013613/1.}

\thanks{$^*$Corresponding author: Christian Bick}

\begin{document}
\maketitle

% Enter the first author's name and email address; email addresses are required for each author.
% Use footnote notations to indicate address and affiliations with commas between numbers if more than one address applies; see below for examples.
\centerline{\scshape
Christian Bick$^{{\href{mailto:c.bick@vu.nl}{\textrm{\Letter}}}*1\text{--}4}$,
Jos\'e Mujica$^{1}$,
Bob Rink$^{1}$},

\medskip

{\footnotesize
% Enter the full affiliation and country name:
% Do not insert commas or periods at the end of lines.
 \centerline{${}^1$Department of Mathematics, Vrije Universiteit Amsterdam, De Boelelaan 1111, Amsterdam, The Netherlands}
} % Do not forget to end {\footnotesize with the sign }

\medskip

{\footnotesize
 % Enter the full affiliation and country name:
 \centerline{$^{2}$Institute for Advanced Study, Technical University of Munich, Lichtenbergstr 2, 85748 Garching, Germany}
}

\medskip

{\footnotesize
 % Enter the full affiliation and country name:
 \centerline{${}^3$Department of Mathematics, University of Exeter, Exeter EX4 4QF, United Kingdom}
}

\medskip

{\footnotesize
 % Enter the full affiliation and country name:
 \centerline{${}^4$Mathematical Institute, University of Oxford, Oxford OX2 6GG, United Kingdom}
}

%\bigskip

% The name of the handling editor will be entered by AIMS production staff.
% "Communicated by Handling Editor" is not needed for special issue.
%\centerline{(Communicated by Handling Editor)}

%%%%%%%%%%%%%%%%%%%%%%%%%%%%%%%%%%%%%%%%%%%%%%%%%%%%%%%
%             5. ABSTRACT
%%%%%%%%%%%%%%%%%%%%%%%%%%%%%%%%%%%%%%%%%%%%%%%%%%%%%%%

\begin{abstract}
We consider synchrony patterns in coupled phase oscillator networks that correspond to invariant tori. 
For specific nongeneric coupling, these tori are equilibria relative to a continuous symmetry action.
We analyze how the invariant tori deform under forced symmetry breaking as more general network interaction terms are introduced.
We first show in general that perturbed tori that are relative equilibria can be computed using a parametrization method; this yields an asymptotic expansion of an embedding of the perturbed torus, as well as the local dynamics on the torus. 
We then apply this result to a coupled oscillator network, and we numerically study the dynamics on the persisting tori in the network by looking for bifurcations of their periodic orbits in a boundary-value-problem setup. 
This way we find new bifurcating stable synchrony patterns that can be the building blocks of larger global structures such as heteroclinic cycles.
\end{abstract}

%%%%%%%%%%%%%%%%%%%%%%%%%%%%%%%%%%%%%%%%%%%%%%%%%%%%%%
%                   6. BODY
%%%%%%%%%%%%%%%%%%%%%%%%%%%%%%%%%%%%%%%%%%%%%%%%%%%%%%

% Only the first word and proper nouns of section titles should be capitalized.
% The title of section 1:

\section{Introduction}

Many real-world physical systems, including biological and neural systems, can be described in terms of interacting oscillatory processes (see, for example,~\cite{Breakspear17,Strogatz2004}).  
Networks of coupled oscillators---in which each node in the network represents an oscillator and the network structure organizes the interaction among the nodes---provide mathematical models for such systems. The emergence of collective behavior, such as (partial) synchronization, is a remarkable effect of the interaction of oscillatory units in coupled oscillator networks; cf.~\cite{BB68,Pikovsky2003,SAMEO05}.

Phase oscillator networks are an important class of oscillator networks in which each oscillator is described by a single phase variable~\cite{ABPRS06,RPJK16,Strogatz00}.
More specifically, the phase~$\theta_k$ of oscillator~$k$ lies on a torus~$\T:=\mathbb R/2\pi\mathbb Z$ and for a network of $n$~oscillators, the joint state $\theta = (\theta_1, \dotsc, \theta_n)$ thus lies on an~$n$-torus~$\T^n$. 
Full (phase) synchrony corresponds to the set~$\Sp := \{\theta_1 = \dotsb =\theta_n\}$, where all oscillators have the same phase.
By contrast, the splay state~$\Dp := \{\theta_k=\theta_{k-1}+\frac{2\pi}{n}\}$ corresponds to a desynchronized configuration where all oscillators are distributed uniformly across the circle~$\T$.
These phase configurations are dynamically invariant for all-to-all coupled networks of identical (generalized) Kuramoto oscillators, that evolve according to
\begin{equation}
\dot\theta_k = \omega + \sum_{j=1}^n g(\theta_j-\theta_k).
\label{eq:singlepop}
\end{equation}
Here~$\omega\in\R$ is the intrinsic frequency of the oscillators and the~$2\pi$-periodic function~$g$ determines pairwise phase interaction between them.
Specifically, the sets~$\Sp$ and~$\Dp$ are equilibria relative to the symmetry action $\theta_j\mapsto \theta_j+\phi$ that shifts the phases of all oscillators by a common constant~$\phi\in\T$.

Many relevant networks are not all-to-all coupled but are rather organized into distinct communities or populations.
In a network consisting of $m$~populations of $n$~oscillators each, let~$\theta_{\sigma,k}\in\T$ denote the phase of the~$k$th oscillator of population~$\sigma\in\{1,\dotsc, m\}$.
Denoting by $\theta_\sigma\in\T^n$ the joint state of population~$\sigma$, then each population may be synchronized or in splay state.
Hence, there may be localized synchrony patterns (or coherence-incoherence patters), where some populations are phase synchronized while others are desynchronized; these patterns are sometimes also referred to as chimeras~\cite{O18,Panaggio2015}.
For instance, for an oscillator network consisting of $m=3$ populations the first population may be synchronized, $\theta_1\in\Sp$, while the others are in splay state, $\theta_2,\theta_3\in\Dp$; we denote this configuration by~$\SDD$ and extend the notation to other synchrony patterns such as~$\SSD$.
Since these synchronization patterns only specify the phase relationship within populations, they are $m$-dimensional tori.

These synchrony patterns arise naturally as invariant sets in coupled populations of phase oscillators, where~$\theta_{\sigma,k}$ evolves according to
\begin{equation}
\dot\theta_{\sigma,k} = \omega + \displaystyle\sum_{j=1}^n g(\theta_{\sigma,j}-\theta_{\sigma,k}) + Y_{\sigma,k}(\theta) =: H_{\sigma,k}(\theta)
\label{eq:multipop}
\end{equation}
where~$g$ mediates the coupling within populations, while the functions $Y(\theta)=\{Y_{\sigma,k}(\theta)\}$ determine the coupling between populations.
We assume in this paper that the coupling function~$Y_{\sigma,k}$ only depends on within-population phase differences, i.e., variables of the form
$\theta_{\tau,j}-\theta_{\tau,k}$.
This ensures that~\eqref{eq:multipop} is equivariant with respect to a $\T^m$-symmetry action, where $\phi = (\phi_1, \dotsc, \phi_m)$ acts by
\begin{align}
\label{eq:Aaction}
\theta_{\sigma, k} \mapsto \theta_{\sigma, k} + \phi_{\sigma}\, .
\end{align}
In other words, there are $m$~independent phase-shift symmetries, each of which shifts the phases  of one particular population by the same constant---analogous to the single phase shift symmetry for one population~\eqref{eq:singlepop}. 
Synchronization patterns, such as~$\SDD$ and~$\SSD$ for $m=3$ populations, are equilibria relative to the phase-shift symmetry action~\eqref{eq:Aaction}.

These synchrony patterns can be part of larger structures that shape the global network dynamics.
For instance, in~\cite{Bick2019} the authors show the existence of heteroclinic cycles between synchrony patterns that are relative equilibria in a system of the form~\eqref{eq:multipop}; the saddle sets involved in the heteroclinic structure are invariant tori of the form~$\SDD$,~$\SSD$, etc.
The key point here is that the independent phase-shift symmetries make the inter-population coupling functions~$Y_{\sigma,k}$ \emph{nongeneric}.
More generic couplings will break these symmetries; this is also called forced symmetry breaking~\cite{Guyard1999,Guyard2001}. In this paper we will consider forced symmetry breaking in systems of the form
\begin{equation}
\begin{split}
\dot\theta_{\sigma,k} &= \omega + \displaystyle\sum_{j=1}^n g(\theta_{\sigma,j}-\theta_{\sigma,k}) + Y_{\sigma,k}(\theta)+\delta Z_{\sigma,k}(\theta)\\%[0.5cm]
&= H_{\sigma,k}(\theta)+\delta Z_{\sigma,k}(\theta)\, ,
\end{split}
\label{eq:multipop_pert}
\end{equation}
where~$Z=\{Z_{\sigma,k}\}$ is a perturbation that breaks the~$\T^m$ symmetry of system~\eqref{eq:multipop} for~$\delta>0$. 
A natural choice for~$Z$ is a function depending on terms of the form~$\theta_{\tau,j}-\theta_{\sigma,k}$, i.e., a phase difference coupling as in~\eqref{eq:singlepop} between \emph{any} two oscillators in the network, not just in the same population.
In the example of the three-population network, this breaks the symmetry from~$\T^3$ to~$\T$.

Here we take the first step in understanding how the heteroclinic structures found in~\cite{Bick2019} perturb as symmetry is broken: 
We study the dynamics on the perturbed tori for~\eqref{eq:multipop_pert}.
The main contributions are the following.
First, we develop/generalize perturbation theory based on the parametrization method for relative equilibria. 
These results may also be of general interest beyond systems of the form~\eqref{eq:multipop}.
Second, we employ this theory to compute approximations of the perturbed tori and vector fields on them. 
Third, we use these approximations to numerically study the solutions on these perturbed tori and elucidate their bifurcations. 
To this end we implement a two-point boundary-value-problem setup and continue it with {\sc Auto}~\cite{auto}. 
This gives insight in how trajectories  approach the torus, how they behave near the torus, and where they are expected to leave. 
This can be used further to understand how the actual heteroclinic cycles perturb.

This paper is organized as follows. In \sref{sec:parreleq} we discuss the problem of perturbing a $\T^m$-equivariant dynamical system with a normally hyperbolic relative equilibrium, and show that the parametrization method developed in~\cite{vondergracht2023parametrisation} can be applied in this setting. 
\Sref{sec:network} presents the model equations of a coupled oscillator network  having invariant tori that are relative equilibria; they relate to synchrony patterns in the network and are the saddle invariant sets involved in the heteroclinic cycles described in~\cite{Bick2019}. 
This is the model in which we exploit the ideas from the parametrization method in order to obtain the reduced dynamics on the (embedded) perturbed tori. 
In \sref{sec:Param} we compute the perturbed tori and dynamics thereon by solving the relevant equations from the parametrization method for our coupled oscillator network. We focus on  perturbations that can be written in terms of phase differences and with sine coupling including second harmonics. 
In \sref{sec:dyn_tori} we use the embedding and the normal forms to study the local dynamics on these perturbed tori.
Specifically, we investigate numerically how periodic orbits that foliate the invariant torus in the unperturbed case bifurcate when we perturb the system (and hence the invariant tori). 
We conclude with a brief discussion and outlook. In the appendix we briefly discuss the numerical implementation for the bifurcation analysis of periodic orbits on the perturbed tori.

%%%%%%
\section{Parametrizing perturbed relative equilibria}
\label{sec:parreleq}

In this section, we introduce a method for computing persisting normally hyperbolic relative equilibria in perturbations of equivariant ODEs. The method
 works by iteratively solving a ``conjugacy equation''. Thus, it does not only yield an approximation of the perturbed relative equilibrium, but also of the dynamics on it. We need both for our numerical study in Section \ref{sec:dyn_tori}. 
 
 The method presented here was introduced in \cite{vondergracht2023parametrisation} for non-equivariant problems, i.e., to  compute approximations of general normally hyperbolic quasi-periodic tori. The work in \cite{vondergracht2023parametrisation} was in turn inspired by papers of De la Llave et al.~\cite{CFD03_1,CFD03_2,CFD03_3}, who popularized the idea of computing invariant manifolds by solving a conjugacy equation. 

\subsection{Equivariant dynamics and relative equilibria}  
\label{sec:torper}

Consider a differential equation 
\begin{equation}
\dot \theta = H(\theta)\ \mbox{for}\ \theta \in \T^d\, ,
\label{eq:vf}
\end{equation} 
and assume that $H$ is equivariant under a  symmetry action on $\T^d$ of the form 
\[\theta \mapsto \theta + R \phi \, , \ \mbox{with} \ \theta \in \T^d \ \mbox{and}\ \phi\in \T^m\, .\]
Here $R: \mathbb{R}^m\to \mathbb{R}^d$ is a linear map whose matrix has integer coefficients. This ensures that~$R$ can  be viewed a map from~$\T^m$ to~$\mathbb{T}^d$, so that the action is well-defined.  We also require that the action is free,  which implies that~$R$ is injective.
The vector field $H$ is $\T^m$-equivariant if it satisfies the identity 
\[H(\theta + R \phi) =  H(\theta)\ \mbox{for all}\ \theta \in \mathbb{T}^d\ \mbox{and all} \ \phi\in \mathbb{T}^m  \, .\]
This equivariance implies that whenever~$\theta(t)$ is a solution of~\eqref{eq:vf}, then so is the curve $t\mapsto \theta(t) + R\phi$, for any $\phi\in \T^m$. 
In other words, solutions to \eqref{eq:vf} come in $\T^m$-orbits.
We refer to~\cite{Golubitsky} for more details on equivariant dynamical systems. 

Our motivating example consists of the coupled populations of Kuramoto oscillators, with~$H$ defined as in~\eqref{eq:multipop}, $d=m\times n$, and the map~$R$ as given in formula \eqref{eq:Aaction}, that is,
\begin{equation}
\label{eq:Rot}
    R \phi = (\phi_1, \dotsc, \phi_1; \phi_2, \dotsc, \phi_2; \dotsc; \phi_m, \dotsc, \phi_m)\in (\T^n)^m\, .
\end{equation}
Next, we assume that $H$ possesses a relative equilibrium: a group orbit that is invariant under the flow of $H$. In the example of the coupled populations of Kuramoto oscillators, a group orbit is a set of configurations with constant within-population phase difference, e.g., the coherence-incoherence states in which some populations are synchronized while others are in a splay state. 

In general, $\theta_0 \in \T^d$ lies on a relative equilibrium when there is
a vector $\Omega \in \mathbb{R}^m$ such that
\[H(\theta_0) = R \Omega \in \, {\rm im}\, R = {\bf T}_{\theta_0}(\T^m\cdot \theta_0)\, .\]
That is, $H(\theta_0)$ lies in the tangent space to the group orbit of $\theta_0$. The $\T^m$-equivariance  implies that 
$H(\theta_0 + R\phi) = H(\theta_0) = R\Omega$ as well, and therefore~$H$ is tangent to the group orbit at~$\theta_0$, if and only if it is tangent to this group orbit everywhere. The group orbit is then an invariant manifold for~$H$. 
The map
\begin{align}\label{def:E}
E: \T^m\hookrightarrow\T^d \ \mbox{defined by}\ 
 E(\phi) := \theta_0 + R\phi 
 \end{align}
 is an explicit embedding of this invariant manifold (recall that the action is free so~$E$ is injective). 
 In fact, we have
 \begin{proposition}\label{prop:Embed}
    The embedding~$E$ defined in~\eqref{def:E} semi-conjugates the constant vector field~$\Omega$ to~$H$, i.e., $E$~sends solutions of $
\dot \phi=\Omega  \ \mbox{on}\ \T^m
$
to solutions of~\eqref{eq:vf}.
 \end{proposition}
 \begin{proof}[Proof of Proposition~\ref{prop:Embed}]
    Let $\phi(t)\in\T^m$ satisfy $\dot \phi = \Omega$. Then 
    \[
    \frac{d}{dt} E(\phi(t)) =R \dot \phi(t) = R \Omega = H(\theta_0)=H(\theta_0 + R\phi(t)) = H(E(\phi(t)))\, .
    \]
    So the result is an immediate consequence of the equivariance of~$H$.
 \end{proof} 
 \noindent Note that we use the word ``semi-conjugacy'' as it is used in \cite{parameterization_book}: a non-invertible map sending solutions to solutions (in this case an embedding).
 
 Depending on the resonance properties of $\Omega$, the solutions of $\dot \phi = \Omega$ are periodic or quasi-periodic  curves: $\phi(t) = \phi_0 + \Omega t \in \mathbb{T}^m$. We conclude that a relative equilibrium of a $\T^m$-action is always a (quasi-)periodic torus. (This is actually a general fact that is not only true in the setting presented here.)

\subsection{A symmetry breaking perturbation}

\noindent Next, we consider a perturbation of~\eqref{eq:vf} of the form
\begin{equation}
\dot{\theta} = F(\theta) := H(\theta)+\delta Z(\theta)\ \ \mbox{with}\  0< \delta \ll 1\, .
\label{eq:pert_gen}
\end{equation}
 As before, we assume that $H$ is $\T^m$-equivariant and that it possesses a relative equilibrium with an embedding $E:\T^m\to\T^d$ of the form given in \eqref{def:E}. 
We do not assume that~$Z$ is $\T^m$-equivariant, i.e., the perturbation  breaks the symmetry of~$H$. 
When $E(\T^m)$ is normally hyperbolic for the unperturbed dynamics~\eqref{eq:vf}, then the relative equilibrium will persist in the perturbation~\eqref{eq:pert_gen} as an invariant torus by F\'enichel's theorem~\cite{Fen71}. 
As in~\cite{vondergracht2023parametrisation}, we compute this persisting torus by searching for a perturbed embedding
\[
e = E +\delta e_1 + \delta^2 e_2 + \mathcal{O}(\delta^3) :\T^m\hookrightarrow\T^d
\]
and a perturbed ``reduced'' vector field 
\[
f =\Omega + \delta f_1 + \delta^2 f_2 + \mathcal{O}(\delta^3) \ :\T^m\to \R^m\, ,
\]
which we require  to satisfy the conjugacy equation
\begin{equation}
De(\phi) f(\phi)=F(e(\phi))\, .
\label{eq:conj_per}
\end{equation}
A solution~$(e,f)$ to~\eqref{eq:conj_per} consists of an embedding~$e$ of an invariant torus for~$F$, and a reduced vector field~$f$ on~$\T^m$ that is semi-conjugate to~$F$. 
This means that~$e$ sends solutions to $\dot \phi=f(\phi)$ to solutions to $\dot \theta =F(\theta)$ and hence~$e(\T^m)$ is invariant for the dynamics of~$F$. 

As shown in~\cite{vondergracht2023parametrisation}, expanding equation~\eqref{eq:conj_per} in~$\delta$ yields a list of iterative equations for $e_1, f_1, e_2, f_2$, etc. The first of these equations is simply $R \Omega = H(E(\phi))$, which holds by assumption. 
The second equation---that is, the $\mathcal{O}(\delta)$-part of \eqref{eq:conj_per}---reads
 $R f_1(\phi) + De_1(\phi)  \Omega = DH(E(\phi)) e_1(\phi) + Z(E(\phi))$. We rewrite this as
\begin{equation} 
(\partial_{\Omega}-DH\circ E) e_1+R  f_1=Z\circ E =:G_1\ .
\label{eq:main_eq}
\end{equation} 
We view~\eqref{eq:main_eq} as a linear inhomogeneous equation for the pair $(e_1, f_1)$. The inhomogeneous right hand side~$G_1$ is a map from~$\T^m$ to~$\R^d$. 
The operator~$\partial_{\Omega}$ denotes  differentiation in the direction of the vector $\Omega$, i.e., $(\partial_{\Omega} e)(\phi) := De(\phi)  \Omega$. 

We refer to~\eqref{eq:main_eq} as the {\it infinitesimal conjugacy equation}, as it is the linear approximation of~\eqref{eq:conj_per}. Higher-order equations (for $e_2, f_2$, etc.) can similarly be derived, see~\cite{vondergracht2023parametrisation}, but we will not use them in this paper. One of our goals  will be to solve equation~\eqref{eq:main_eq} for some examples, to obtain information on the location and the dynamics of the perturbed torus.  

\subsection{The linearized dynamics around a relative equilibrium}
In order to solve equation~\eqref{eq:main_eq}, we will decompose the perturbation $e_1$ of the embedding~$e$ into a component tangent to the unperturbed torus~$E(\T^m)$ and a component normal to it. To do this in a meaningful way, we need a better understanding of the linearized dynamics near this unperturbed torus. 
 In fact, it was shown in~\cite{vondergracht2023parametrisation} that it is crucial that the unperturbed torus is {\it reducible}. We shall not introduce this concept here, because we can prove something slightly stronger for relative equilibria, as a result of the following simple observation.
\begin{proposition}\label{lemmarelativeequilibrium}
Let the vector field~$H$ on~$\T^d$ be $\T^m$-equivariant as above. Then the Jacobian matrix~$DH(\theta)$ is constant along group orbits, and in particular along relative equilibria, i.e., $DH \circ E = DH(\theta_0)$.
\end{proposition}
\begin{proof}[Proof of Proposition~\ref{lemmarelativeequilibrium}]
    Differentiating $H(\theta + R \phi) = H(\theta)$ with respect to~$\theta$ gives $DH(\theta+R\phi) = DH(\theta)$, showing that $DH$ is constant along any group orbit. Setting $\theta = \theta_0$ gives $DH(E(\phi))=DH(\theta_0)$.
\end{proof}
\noindent The next result characterizes the linearized dynamics of~$H$ tangent to the relative equilibrium.
\begin{proposition}\label{prop:tangentialdynamics}
   We have $DH(\theta) R = 0$ for all $\theta\in \T^d$. In other words, 
   $${\rm im}\, R \subset \ker DH(\theta)\, .$$  
\end{proposition}
\begin{proof}[Proof of Proposition~\ref{prop:tangentialdynamics}]
    Differentiating the equivariance identity $H(\theta + R\phi) =  H(\theta)$ with respect to~$\phi$ at $\phi=0$ gives the result.  
\end{proof}
\noindent Applied to $\theta = \theta_0+R\phi \in E(\T^m)$, Proposition~\ref{prop:tangentialdynamics} reduces to the statement that the tangent space ${\rm im}\, R$ to the relative equilibrium is in the kernel of $DH$, i.e.,  \[(DH \circ E)R=0\, .\]
We  need a similar simple description of the  linearized  dynamics normal to the relative equilibrium. In case the relative equilibrium is normally hyperbolic, this description is given by the following proposition. 
\begin{proposition}\label{prop:splittingprop}
    Assume that the relative equilibrium $E(\T^m)\subset \T^d$ is normally hyperbolic for the unperturbed dynamics \eqref{eq:vf}. Then there exist an injective linear map $N:  \R^{d-m}\to \R^d$ and a hyperbolic linear map $L: \mathbb{R}^{d-m} \to \mathbb{R}^{d-m}$  such that 
      \begin{itemize}
          \item[{\bf 1.}]
          The range of~$N$ is normal to the relative equilibrium:
      $$\R^d = {\rm im}\, R \oplus {\rm im}\, N\, .$$
      \item[{\bf 2.}] $N$ semi-conjugates $L$ to $DH$ on the relative equilibrium:
      $$(DH\circ E) N = N L\, .$$
      \end{itemize}
\end{proposition}

\begin{proof}[Proof of Proposition~\ref{prop:splittingprop}]
Recall that ${\rm im}\, R \subset  \ker DH \circ E = \ker DH(\theta_0)$. 
    The assumption that the relative equilibrium is normally hyperbolic implies that $\ker DH(\theta_0) = {\rm im}\, R$, and also  that this kernel can be complemented by the sum of the hyperbolic generalized eigenspaces of $DH(\theta_0)$. Let $N: \R^{d-m} \to\R^d$ be an injective linear map spanning the sum of these generalized eigenspaces (for instance, the matrix~$N$ could have the hyperbolic generalized eigenvectors as its columns). By construction it then holds that $\R^d = {\rm im}\, R \oplus {\rm im}\, N$.
   It also follows that there exists a linear map $L: \R^{d-m} \to\R^{d-m}$ satisfying 
    $$DH(\theta_0) N = N L\, .$$  
    The eigenvalues of~$L$ are the nonzero eigenvalues of $DH(\theta_0)$ so $L$~is hyperbolic. Moreover, by Proposition~\ref{lemmarelativeequilibrium} we have $(DH \circ E)N = DH(\theta_0)N=NL$.
\end{proof}
\noindent The injective matrix $N$ in Proposition \ref{prop:splittingprop} is clearly not unique: we only require that its image is the hyperbolic subspace of $DH(\theta_0)$. As a result, $L$ is not unique either, although the eigenvalues of $L$ are fixed and equal to the nonzero eigenvalues of $DH(\theta_0)$. A matrix $L$ satisfying the conclusion of the proposition is called a \emph{Floquet matrix} for the invariant torus $E(\mathbb{T}^m)$. Its (unique) eigenvalues are called the \emph{Floquet exponents} of the invariant torus. The normal hyperbolicity of the torus implies that none of these Floquet exponents lie on the imaginary axis.

\subsection{A solution to the infinitesimal conjugacy equation}
\label{sec:ConjSol}
Proposition~\ref{prop:splittingprop} enables us to solve the infinitesimal conjugacy equation~\eqref{eq:main_eq} by decomposing the perturbed embedding into a component tangent to and a component normal to the relative equilibrium. 
Specifically, we make the ansatz
\begin{align}\label{eq:cleveransatz}
e_1(\phi) = R X_1(\phi) + N Y_1(\phi)\, .
\end{align}
Here, the functions
\[X_1:\T^m \to \R^m \ \mbox{and} \ Y_1:\T^m\to \R^{d-m}\,\]
are still to be determined. The following fundamental proposition can also be found in~\cite{vondergracht2023parametrisation}.
\begin{proposition}\label{prop:Ansatz}
    The ansatz~\eqref{eq:cleveransatz} transforms equation~\eqref{eq:main_eq} into 
    \begin{align}
\label{eq:afteransatz}    
    R ( \partial_{\Omega}X_1 + f_1) + N (\partial_{\Omega} -L) (Y_1) = G_1\, .
    \end{align}
\end{proposition}
\begin{proof}[Proof of Proposition~\ref{prop:Ansatz}]
    The result follows from inserting the ansatz \eqref{eq:cleveransatz} into the infinitesimal conjugacy equation \eqref{eq:main_eq} and using that 
    $(DH\circ E)R=0$, that $(DH \circ E) N=NL$, and that $R$ and $N$ are independent of $\phi\in \T^m$. (We remark that this latter fact is not true in the non-equivariant setting considered in~\cite{vondergracht2023parametrisation}.) 
\end{proof}
\noindent As~$R$ and~$N$ are transverse, equation~\eqref{eq:afteransatz} can be decomposed into two equations by projecting it onto the tangent space and the normal space. To this end, let us denote by $\pi: \mathbb{R}^d \to \mathbb{R}^d$ the projection onto ${\rm im}\, R$ along ${\rm im}\, N$. This means that $\pi$ is the unique linear map satisfying $\pi R = R$ and $\pi N = 0$. Applying~$\pi$ and~$1-\pi$ to~\eqref{eq:afteransatz} yields
\begin{align*}
R( \partial_{\Omega}X_1 + f_1 )&= \pi G_1 \in {\rm im}\, R && \text{and} & N (\partial_{\Omega} -L) (Y_1) &= (1-\pi)G_1 \in {\rm im}\, N\, .     
\end{align*}
Equivalently, we may write these equations as
\begin{subequations}\label{eq:homological}
\begin{align}
 \partial_{\Omega}X_1 + f_1 &= R^+ \pi G_1 =: U_1
 \label{eq:HomlogicalTan}\\   
 (\partial_{\Omega} -L) (Y_1) &= N^+ (1-\pi)G_1 =:V_1\, .       
 \label{eq:HomlogicalNor} 
\end{align}    
\end{subequations}
Here, $R^+ := (R^\trp R)^{-1}R^\trp$ and $N^+ = (N^\trp N)^{-1}N^\trp$ denote the Moore--Penrose pseudo-inverses of~$R$ and~$N$. 

We refer to~\eqref{eq:HomlogicalTan} as the first \emph{tangential homological equation} and to~\eqref{eq:HomlogicalNor} as the first \emph{normal homological equation}.
We see from these homological equations that~$X_1$ can be chosen freely, while~$Y_1$ and~$f_1$ are then given by 
\begin{equation}\label{eq:HomologicalSol}
Y_1 = (\partial_{\Omega} - L)^{-1}V_1 \ \mbox{and}\ f_1 = U_1 - \partial_{\Omega}X_1\, .
\end{equation}
It was shown in~\cite{vondergracht2023parametrisation} that $X_1, Y_1, f_1$ can be found as Fourier series in $\phi\in \T^m$. In particular, the hyperbolicity of $L$ implies that the operator $\partial_{\Omega}-L$ is invertible. It was also shown in  \cite{vondergracht2023parametrisation} that $X_1$ can be chosen in such a way that $f_1$ is in normal form, i.e., that one can remove ``non-resonant terms'' from $f_1$ by choosing $X_1$ appropriately. 
We will not reprove any of these facts here. Instead, we will simply observe them in our computations later on.

%%%%%%
\section{Symmetry breaking in coupled phase oscillator networks}
\label{sec:network}

\newcommand{\K}{K}
\newcommand{\tGf}{G_4}
\newcommand{\gt}{g_2}
\newcommand{\Xp}{\mathrm{X}}

We now exploit these results to elucidate the effect of forced symmetry breaking in $m$~coupled populations of $n$~oscillators given by the perturbation problem~\eqref{eq:multipop_pert}; such phase equations can be motivated by phase reductions~\cite{Ashwin2022,Bick2023}.
The unperturbed equations for $\delta=0$ are equivariant with respect to the action of~$\T^m$ given by~\eqref{eq:Aaction}, which corresponds to a phase-shift symmetry in each population.
This implies that the synchrony patterns where each population is either synchronized ($\theta_\sigma\in\Sp$) or desynchronized ($\theta_\sigma\in\Dp$),
\begin{align}
\Xp_1\dotsc\Xp_m := \{\theta = (\theta_1, \dotsc,\theta_m)\in\T^{mn}:\theta_\sigma\in\Xp_\sigma\in\{\Sp,\Dp\}\}
\label{eq:Tori}
\end{align}
are equilibria relative to the symmetry action of $\T^m$. 
For example, $\Sp\dotsc\Sp$ ($m$-times) is the relative equilibrium where all populations are phase synchronized individually (but not necessarily to each other).
Perturbing the system, $\delta > 0$, breaks the~$\T^m$ symmetry to a~$\T$ symmetry where $\alpha\in\T$ acts by a common phase shift of all phases, $\theta_{\sigma,k}\mapsto \theta_{\sigma,k}+\alpha$.
We can now compute how invariant tori~\eqref{eq:Tori} and the dynamics deform as the symmetry is broken.

While the results apply more generally to systems of the form~\eqref{eq:multipop_pert}, we focus on the case of $m=3$ population of $n=2$ oscillators as considered in~\cite{Bick2019,Bick2017c}.
Specifically, the oscillators in population~$\sigma$ evolve according to
\begin{subequations}\label{eq:Dynamics3x2}
\begin{align}
\dot\theta_{\sigma,1} &= H_{\sigma, 1}(\theta) +\delta Z_{\sigma,1}(\theta) =  \omega + g(\theta_{\sigma,2}-\theta_{\sigma,1}) + Y_{\sigma,1}(\theta)+\delta Z_{\sigma,1}(\theta)\\
\dot\theta_{\sigma,2} &= 
H_{\sigma, 2}(\theta) +\delta Z_{\sigma,2}(\theta)
= \omega + g(\theta_{\sigma,1}-\theta_{\sigma,2}) + Y_{\sigma,2}(\theta)+\delta Z_{\sigma,2}(\theta)
\end{align}
\end{subequations}
where the coupling between populations with strength $K\geq 0$ is determined by
\begin{equation}
\begin{split}
    Y_{\sigma,k}(\theta) &= 
-\frac{\K}{4}
\big(g_4(\theta_{\sigma-1, 1}-\theta_{\sigma-1,2}+\theta_{\sigma,3-k}-\theta_{\sigma,k})
\\&\qquad\qquad
+g_4(\theta_{\sigma-1, 2}-\theta_{\sigma-1,1}+\theta_{\sigma,3-k}-\theta_{\sigma,k})\big)\\
&\ \quad+\frac{\K}{4}
\big(g_4(\theta_{\sigma+1, 1}-\theta_{\sigma+1,2}+\theta_{\sigma,3-k}-\theta_{\sigma,k})
\\&\qquad\qquad
+g_4(\theta_{\sigma+1, 2}-\theta_{\sigma+1,1}+\theta_{\sigma,3-k}-\theta_{\sigma,k})\big)
\end{split}
\end{equation}
for $k=1,2$ and $\sigma=1,2,3$ taken modulo~$3$.
The $2\pi$-periodic coupling functions
\begin{subequations}
\label{eq:coupling}
\begin{align}
g(\vartheta) = g_2(\vartheta)&=\sin(\vartheta+\alpha_2)-r_0\sin(2(\vartheta+\alpha_2)),\\
g_4(\vartheta)&=\sin(\vartheta+\alpha_4),
\end{align}
\end{subequations} with parameters $\alpha_2, \alpha_4 \in \R$ mediate the pairwise interactions within populations and nonpairwise interactions between populations, respectively. 
For $\delta=0$, the system~\eqref{eq:Dynamics3x2} is~$\T^3$ equivariant.

For $\delta>0$ the network dynamics~\eqref{eq:Dynamics3x2} is now subject to forced symmetry breaking.
We focus on a symmetry breaking perturbation of the form
\begin{equation}
Z_{\sigma,k}(\theta)=\sum_{l,m} h(\theta_{l,m}-\theta_{\sigma,k}),
\label{eq:Zj}
\end{equation}
where~$h$ is a (nonconstant) $2\pi$-periodic function.
This corresponds to pairwise Kuramoto-type coupling between any two oscillators as in~\eqref{eq:singlepop} and breaks the phase shift symmetries induced by the nonpairwise coupling between populations: The $\T^3$~symmetry is broken to a single $\T$~symmetry.
Together, Equations~\eqref{eq:Dynamics3x2}--\eqref{eq:Zj} define the system we investigate.

%%%
\subsection{Invariant tori as relative equilibria for $\delta=0$}
\label{sec:RelativeEq}

\newcommand{\Z}{\mathbb{Z}}
\newcommand{\Tbf}{\mathbf{T}}

Because of the phase-shift symmetry, we may assume $\omega=-1$ without loss of generality\footnote{While this value is different from the more common choice of $\omega=0$, it is convenient for plotting as synchronized solutions appear stationary.}. 
In addition to the continuous symmetries, the system also admits discrete symmetries: 
Write $\Z_m := \Z/m\Z$ for the cyclic group of $m$~elements.
The dynamical equations~\eqref{eq:Dynamics3x2} are equivariant with respect to an action of~$\Z_3$ that acts by permuting the population index~$\sigma$ and~$(\Z_2)^3$ that acts by permuting the indices~$k$ within populations.

The unperturbed dynamics for $\delta=0$ have eight invariant tori, $\SSS, \DDD$ as well as $\SSD, \SDD, \DSS, \DSD, \SDS, \DDS$. 
We focus on the latter as they are part of the heteroclinic cycles found in~\cite{Bick2019}.
For concreteness, we set $(\alpha_2,\alpha_4)=(\frac{\pi}{2},\pi)$ in the parameter range where where heteroclinic cycles involving these tori may exist.
We can further restrict our attention to~$\SDD, \SSD$ as the remaining tori are their images under~$\Z_3$ action.

The torus $\SDD$ is a copy of~$\T^3$ embedded into~$\T^6$ through the map
\begin{align}
E^{\tiny \SDD}(\phi_1,\phi_2,\phi_3)&=(\phi_1,\phi_1;\phi_2,\phi_2+\pi;\phi_3,\phi_3+\pi)
\label{eq:e0_SDD}
\end{align}
Since $g_2(0) = 1$, $g_2(\pi) = -1$ and $g_4(0)= g_4(\pi) = 0$,
the constant vector field on $\SDD$ is given by $\Omega^\SDD=(0,-2,-2)^\trp$.
Transverse stability of the invariant torus is governed by the eigenvalues corresponding to the directions that break~$\Sp$, $\Dp$ in each direction (e.g., the eigenvector $(-1, 1; 0, 0; 0, 0)$ for the first population); cf.~\cite{Bick2019}.
Taking derivatives these evaluate to $\lambda_1^\SDD = -4r_0$, $\lambda_2^\SDD = -4r_0+2K$, and $\lambda_3^\SDD = -4r_0-2K$.

Similarly, the torus~$\SSD$ can be parameterized by the embedding
\begin{align}
E^{\tiny \SSD}(\phi_1,\phi_2,\phi_3)&=(\phi_1,\phi_1;\phi_2,\phi_2;\phi_3,\phi_3+\pi).
\label{eq:e0_SSD}
\end{align}
The constant dynamics on the torus are determined by the frequency vector $\Omega^\SSD = (0,0,-2)^\trp$.
Transverse stability is governed by the eigenvalues
$\lambda_1^\SSD = -4r_0+2K$,
$\lambda_2^\SSD = -4r_0-2K$, and
$\lambda_3^\SSD = -4r_0$.

%%%
\subsection{Perturbed tori and dynamics thereon}

Our main goal is to understand how the normally hyperbolic invariant tori $\SSD, \SDD$ for~\eqref{eq:Dynamics3x2} deform under forced symmetry breaking for $\delta>0$.
Thus, for $\Tbf\in\{\SSD,\SDD\}$ we seek (first-order) perturbed embeddings
\begin{align}
e^{\Tbf}_{\delta} &= E^\Tbf+\delta e^{\Tbf}_1 :\T^3\hookrightarrow\T^6
\intertext{and vector field on~$\T^3$ given by}
f^\Tbf_{\delta} &= \Omega^\Tbf + \delta f^\Tbf_1.
\label{eq:f_delta}
\end{align}
To compute the perturbed dynamics and embedding for $\delta>0$, we solve the homological equations~\eqref{eq:homological}. We now state the main results and give details of these computations in the following section.

\begin{proposition}\label{prop:Comp}
Consider the phase oscillator network~\eqref{eq:Dynamics3x2} of $m=3$ populations with $n=2$ oscillators each for $(\alpha_2,\alpha_4)=(\frac{\pi}{2},\pi)$
subject to the perturbation~\eqref{eq:Zj}.

For~$\SDD$ we have $\Omega^\SDD = (0, -2, -2)^\trp$, and a normal form for the reduced dynamics on the corresponding perturbation of the invariant torus~$\SDD$ has the form
\begin{equation}
\begin{split}
\dot \phi_1 &= \delta f_{1;1}^\SDD(\phi_3-\phi_2),\\
\dot \phi_2 &= -2+\delta f_{1;2}^\SDD(\phi_3-\phi_2),\\
\dot \phi_3 &= -2+\delta f^\SDD_{1;3}(\phi_3-\phi_2).
\end{split}
\label{eq:SDD_red_gen}
\end{equation}
where~$f_{1;j}^{\SDD}$ depends on the Fourier coefficients~$h_\ell$ of~$h$ according to~\eqref{eq:f_all_SDD}.

Similarly, for~$\SSD$ we have $\Omega^\SDD = (0, 0, -2)^\trp$ and a normal form reads 
\begin{equation}
\begin{split}
\dot \phi_1 &=\delta f^\SSD_{1;1}(\phi_1-\phi_2),\\
\dot \phi_2 &= \delta f^\SSD_{1;2}(\phi_1-\phi_2),\\
\dot \phi_3 &=-2+\delta f^\SSD_{1;3}(\phi_1-\phi_2).
\end{split}
\label{eq:SSD_red_gen}
\end{equation}
where~$f_{1;j}^{\SSD}$ depends on the Fourier coefficients~$h_\ell$ of~$h$ according to~\eqref{eq:f_all_SSD}.
\end{proposition}

Note that the normal form  depends only on combination angles $\phi_k-\phi_j$ for which the frequencies~$\Omega_k$ and~$\Omega_j$ are resonant. As a result, the  dynamics on either torus is effectively one-dimensional and determined by the dynamics of $\psi:=\phi_3-\phi_2$ (for~$\SDD$) and $\varphi := \phi_2-\phi_1$ (for~$\SSD$). 

A classical case is to consider phase interactions that are truncated Fourier series. 
Computing the functions~$f_1^\Tbf$ explicitly yields the following result.

\begin{proposition}
\label{prop:CompWCoup}
Consider the perturbation problem~\eqref{eq:Dynamics3x2} as above with function
\begin{equation}
h(\varphi) = \sin(\varphi+\alpha)+r\sin(2(\varphi+\beta)).
\label{eq:h2}
\end{equation}

Then the first-order correction to the embedding of~$\SDD$ is given by~\eqref{eq:emb_SDD}
and in terms of~$\phi_1$ and $\psi:=\phi_3-\phi_2$ the dynamics on the perturbed torus is, to first order, 
\begin{equation}
\begin{split}
\dot \phi_1 &= 2\delta( \sin\alpha+r\sin(2\beta)),\\
\dot \psi &= 4r\delta(\sin (2(-\psi+\beta))-\sin(2(\psi+\beta))).
\end{split}
\label{eq:SDD_red_lem}
\end{equation}
Similarly, for~$\SSD$ the first-order correction of the embedding of is given by~\eqref{eq:emb_SSD}
and in terms of $\varphi:=\phi_2-\phi_1$ and~$\phi_3$ the dynamics on the perturbed torus is, to first order, 
\begin{equation}
\begin{split}
\dot \varphi &= -4\delta\sin\varphi(\cos \alpha+2r\cos\varphi\cos(2\beta)),\\
\dot \phi_3 &= -2+2r\delta\sin(2\beta).
\end{split}
\label{eq:SSD_red_lem}
\end{equation}
\end{proposition}

These dynamical equations, together with the approximations of the embeddings of SSD and SDD, will allow us to generate starting data to explore solutions of the symmetry broken system numerically. We will do this in Section~\ref{sec:dyn_tori} below.

%%%%%%%
\section{Computing approximations for perturbed tori}
\label{sec:Param}

We now give details of the computations that lead to Propositions~\ref{prop:Comp} and~\ref{prop:CompWCoup}. 
Using the methodology outlined in Section~\ref{sec:parreleq}, we will solve the first order conjugacy equation \eqref{eq:main_eq}, not only for the approximate reduced dynamics $f_1$, but also for the first order approximation of the embedding $e_1$, as we need both  for the numerical continuation in the following section.

%%%
\subsection{Preliminaries}

To compute the first-order approximation of the perturbed tori, we solve the homological equations~\eqref{eq:homological}. We first collect the basic ingredients from Section~\ref{sec:RelativeEq}.
Written as a matrix, the map~$R$ in~\eqref{eq:Rot} that determines the phase-shift symmetries of the unperturbed system---and thus directions tangent to the unperturbed torus---is
\begin{equation}
R = \left(
\begin{matrix}
1&0&0\\
1&0&0\\
0&1&0\\
0&1&0\\
0&0&1\\
0&0&1
\end{matrix}
\right).
\label{eq:De0}
\end{equation}
Because the Jacobian matrices of \eqref{eq:Dynamics3x2} (for $\delta=0$) at SSD and SDD are symmetric matrices, the images of $R$ and $N$ are orthogonal. An example of a matrix $N$ of which the image is orthogonal to that of $R$ is given by
\begin{equation}
N=\left(
\begin{matrix}
-1&0&0\\
1&0&0\\
0&-1&0\\
0&1&0\\
0&0&-1\\
0&0&1
\end{matrix}
\right).
\label{eq:N0}
\end{equation}
The columns of~$N$ represent a ``splitting'' of each oscillator population. The Floquet matrix that determines the linear dynamics in the normal directions is now computed by solving the linear equation $(DH) N = N L$, where~$DH$ is the Jabobian of the unperturbed ODE  given in \eqref{eq:Dynamics3x2}, evaluated at~$\SSD$ or~$\SDD$. 

The entries of the Floquet matrices relate directly to the transverse stability of the relative equilibria (given in Section~\ref{sec:RelativeEq}).
Thus, for the unperturbed torus~$\SDD$, we compute
\begin{equation}
L^\SDD=\left(
\begin{matrix}
-4r_0&0&0\\
0&-4r_0+2K&0\\
0&0&-4r_0-2K
\end{matrix}
\right).
\label{eq:M0_SDD}
\end{equation}
Since the unperturbed dynamics are given by the frequency vector $\Omega^\SDD=(0,-2,-2)^\trp$, populations two and three are resonant. 
Thus, a normal form for the reduced dynamics on the perturbation only depends on~$\phi_1$ and the phase difference~$\phi_3-\phi_2$; cf.~\eqref{eq:SDD_red_gen} in Proposition~\ref{prop:Comp}.
Similarly, for~$\SSD$ we obtain in a similar way the Floquet matrix
\begin{equation}
L^\SSD=\left(
\begin{matrix}
-4r_0+2K&0&0\\
0&-4r_0-2K&0\\
0&0&-4r_0
\end{matrix}.
\right)
\label{eq:M0_SSD}
\end{equation}
With the frequency vector $\Omega^\SSD = (0,0,-2)^\trp$ the normal form equations for the reduced dynamics will only depend on~$\phi_2-\phi_1$ and~$\phi_3$.

\newcommand{\sm}{\smallsetminus}

%%%
\subsection{Perturbing~$\SDD$}
\label{sec:embedding_SDD}

We first focus on the first-order perturbation of~$\SDD$ for the dynamics~\eqref{eq:Dynamics3x2} with perturbation~\eqref{eq:Zj} with global phase-difference coupling determined by the $2\pi$-periodic function~$h$. 
With $i:=\sqrt{-1}$ write $h(\varphi)=\sum_{\ell\in\Z}h_\ell e^{i\ell\varphi}$ with Fourier coefficients $h_\ell \in\mathbb{C}$.
To simplify notation, we will typically suppress the indices ``$\SDD$'' (the torus) and ``$1$'' (the order of the approximation in terms of the small parameter) during the computations.

We first compute the normal component~$Y$ of the embedding. Expanding into a Fourier series and evaluating the right hand side~$V$ of~\eqref{eq:HomlogicalNor} yields
\begin{equation}
V(\phi)=\sum_{l\in\mathbb Z}\left(V_\ell^{12}e^{i\ell(\phi_1-\phi_2)}+V_\ell^{13}e^{i\ell(\phi_1-\phi_3)}\right)
\label{eq:V_SDD}
\end{equation}
with
\begin{align}
V_\ell^{12}&=\left(
\begin{matrix}
0\\
-4h_\ell\\
0
\end{matrix}
\right), &
V_\ell^{13}&=\left(
\begin{matrix}
0\\
0\\
-4h_\ell
\end{matrix}
\right).
\label{eq:Vcoeff}
\end{align}
for $\ell$~odd and
and $V_\ell^{12} = 0$ for $\ell$~even.
With~\eqref{eq:HomologicalSol} this yields
\begin{equation}\label{eq:Y}
Y(\phi)=\sum_{\ell\in\mathbb Z} \left(Y_\ell^{12}e^{i\ell(\phi_1-\phi_2)}+Y_\ell^{13}e^{i\ell(\phi_1-\phi_3)}\right)
\end{equation}
with coefficients
\begin{align}\label{eq:Ycoeff}
Y_\ell^{12}&=\left(
\begin{matrix}
0\\
\dfrac{-2(2r_0-K-i\ell)}{(2r_0-K)^2+\ell^2}h_\ell\\
0
\end{matrix}
\right), &
Y_\ell^{13}=\left(
\begin{matrix}
0\\
0\\
\dfrac{-2(2r_0+K-i\ell)}{(2r_0+K)^2+\ell^2}h_\ell
\end{matrix}
\right).
\end{align}

Next, we consider~\eqref{eq:HomlogicalTan} to obtain the tangential component. With $G=Z\circ E$, the right hand side evaluates to
\begin{equation}
U(\phi)=\sum_{\ell\in\mathbb Z} \left(U_\ell^{00} + U_\ell^{21}e^{i\ell(\phi_2-\phi_1)}+U_\ell^{31}e^{i\ell(\phi_3-\phi_1)}+U_\ell^{32}e^{i\ell(\phi_3-\phi_2)}\right),
\label{eq:U_sdd}
\end{equation}
where
\begin{align}\label{eq:Ul_sdd}
U_\ell^{00} &=\left(
\begin{matrix}
2h_\ell\\
2h_\ell\\
2h_\ell
\end{matrix}
\right), &
U_\ell^{21} &=\left(
\begin{matrix}
2h_\ell\\
2h_{-\ell}\\
0
\end{matrix}
\right),&%\\
U_\ell^{31}&=\left(
\begin{matrix}
2h_\ell\\
0\\
2h_{-\ell}
\end{matrix}
\right), &
U_\ell^{32}&=\left(
\begin{matrix}
0\\
2h_\ell\\
2h_{-\ell}
\end{matrix}
\right)
\end{align}
for $\ell$~even and $U_\ell^{00} =\left(
2h_\ell,
0,
0\right)^\trp$, $U_\ell^{21}=U_\ell^{31}=U_\ell^{32}=0$ for $\ell$~odd.
We can now choose the tangential component of the embedding~$X$ to contain the nonresonant terms such that the first order approximation~$f$ is in normal form with only resonant terms; cf.~\Sref{sec:ConjSol}. 
Specifically, if $\langle \Omega, \lb \rangle = \Omega_1 \ell_1 + \Omega_2 \ell_2 + \Omega_3 \ell_3$ then
\begin{equation}
X(\phi)
=\underset{\tiny{\langle\Omega,\lb \rangle\neq0}}{\sum_{\lb \in\mathbb Z^3}}\frac{1}{i\langle \Omega,\lb  \rangle}U_\lb e^{i\langle \lb ,\phi\rangle },
\label{eq:X}
\end{equation}
only contains nonresonant terms for which $\langle\Omega,\lb \rangle\neq0$.
With this choice of~$X$ and~\eqref{eq:HomologicalSol} we obtain the first order correction of the vector field on the perturbed torus
\begin{equation}
f_1^\SDD(\phi)=\sum_{\ell\in\mathbb Z} \left(U_\ell^{00} + U_\ell^{32}e^{i\ell(\phi_3-\phi_2)}\right)
\label{eq:f_all_SDD}
\end{equation}
in normal form with coefficients as in~\eqref{eq:Ul_sdd}. 
Summarizing, we computed that 
$e^{\SDD}_{\delta} = E^\SDD+\delta e^{\SDD}_1 + \mathcal{O}(\delta^2)$ where 
$e^{\SDD}_1(\phi) = RX(\phi)+NY(\phi)$ for $X(\phi)$ and~$Y(\phi)$ as given in  \eqref{eq:X}, \eqref{eq:Y}, \eqref{eq:Ycoeff}.
For the dynamics on this torus we found that $f^\SDD_{\delta} = \Omega^\SDD + \delta f_1^\SDD + \mathcal{O}(\delta^2)$ with $f_1^\SDD$ as given in \eqref{eq:f_all_SDD}.
This proves the claims of Proposition~\ref{prop:Comp} for~$\SDD$.

We can write out the first-order approximation of the embedding and the dynamics on the perturbed torus for an explicit coupling function.
Since the nontrivial Fourier coefficients $U_\ell^{21}, U_\ell^{31}=U_\ell^{32}$ vanish for odd~$\ell$ of the perturbation function~$h$, we expect that only the second harmonic appears in the reduced dynamics for the biharmonic coupling function
$h(\varphi) = \sin(\varphi+\alpha)+r\sin(2(\varphi+\beta))$
in Proposition~\ref{prop:CompWCoup}.
Indeed, for the coupling function we have
\begin{align*}
h_{-1}&=-\frac{e^{-i\alpha}}{2i}, &
h_{1}&=\frac{e^{i\alpha}}{2i},&
h_{-2}&=-r\frac{e^{-2i\beta}}{2i}, &
h_{2}&=r\frac{e^{2i\beta}}{2i}.
\end{align*}
Evaluating~\eqref{eq:f_all_SDD} explicitly yields phase dynamics
\begin{subequations}\label{eq:nf_SDD}
\begin{align}
\dot\phi_1&= 2\delta(\sin\alpha+r\sin(2\beta))\\
\dot\phi_2&=-2+2r\delta(\sin(2\beta)+\sin(2(\phi_3-\phi_2+\beta))-\sin(2(\phi_2-\phi_3+\beta)))\\
\dot\phi_3&=-2+2r\delta(\sin(2\beta)-\sin(2(\phi_3-\phi_2+\beta))+\sin(2(\phi_2-\phi_3+\beta))).
\end{align}
\end{subequations}
These are embedded in~$\T^6$ via $e^{\SDD}_{\delta} = E^\SDD+\delta e^{\SDD}_1$
with
\begin{equation}\label{eq:emb_SDD}
    e^\SDD_1(\phi)=
\left(
\begin{matrix}
 \frac12r(\cos(2(\phi_2-\phi_1+\beta))+\cos(2(\phi_3-\phi_1-\beta))) \\[0.2cm]
 \frac12r(\cos(2(\phi_2-\phi_1+\beta))+\cos(2(\phi_3-\phi_1-\beta))) \\[0.2cm]
-\frac12r\cos(2(\phi_2-\phi_1+\beta))-a\cos(\phi_1-\phi_2)-b\sin(\phi_1-\phi_2)\\[0.2cm]
-\frac12r\cos(2(\phi_2-\phi_1+\beta))+a\cos(\phi_1-\phi_2)+b\sin(\phi_1-\phi_2)\\[0.2cm]
  -\frac12r\cos(2(\phi_3-\phi_1+\beta))-c\cos(\phi_1-\phi_3)-d\sin(\phi_1-\phi_3)\\[0.2cm]
  -\frac12r\cos(2(\phi_3-\phi_1+\beta))-c\cos(\phi_1-\phi_3)-d\sin(\phi_1-\phi_3)
\end{matrix}
\right)
\end{equation}
and coefficients
\begin{align*}\label{eq:const_Y_SDD}
a&=\dfrac{2(\cos \alpha -(2r_0-K)\sin\alpha)}{(2r_0-K)^2+1}, &
b&=\dfrac{-2(\sin \alpha +(2r_0-K)\cos\alpha)}{(2r_0-K)^2+1}, \\
c&=\dfrac{2(\cos \alpha -(2r_0+K)\sin\alpha)}{(2r_0+K)^2+1}, &
d&=\dfrac{2(\sin \alpha +(2r_0+K)\cos\alpha)}{(2r_0+K)^2+1}.
\end{align*}
This is the main claim of Proposition~\ref{prop:CompWCoup} for~$\SDD$ with $\psi = \phi_3-\phi_2$.

%%%
\subsection{Perturbing~$\SSD$}

An analogous calculation can be done for~$\SSD$; we just state the main results of the computations. For the normal direction, we compute
\begin{equation}
V(\phi)=\sum_{\ell\in\mathbb Z}\left(V_\ell^{13}e^{i\ell(\phi_1-\phi_3)}+V_\ell^{23}e^{i\ell(\phi_2-\phi_3)}\right)
\label{eq:V_SSD}
\end{equation}
with
\begin{align*}
V_\ell^{13} =
V_\ell^{23}&=\left(
\begin{matrix}
0\\
0\\
-4h_\ell
\end{matrix}
\right).
\end{align*}
for $\ell$~odd and $V_\ell^{13}=V_\ell^{23}=0$ for $\ell$~even.
This gives the normal component of the embedding,
\begin{equation}
Y(\phi)=\sum_{\ell\in\mathbb Z}\left(Y_\ell^{13}e^{i\ell(\phi_1-\phi_3)}+Y_\ell^{23}e^{i\ell(\phi_2-\phi_3)}\right)
\end{equation}
with Fourier coefficients
\begin{equation}
Y_\ell^{13}=Y_\ell^{23}=\left(
\begin{matrix}
0\\
0\\
\left(\dfrac{-4(2r_0-i\ell)}{4+\ell^2}\right) h_\ell
\end{matrix}
\right).
\end{equation}
for $\ell$~odd and $Y_\ell^{13}=Y_\ell^{23}=0$ for $\ell$~even.

For the tangential direction we compute the Fourier coefficients
\begin{subequations}
\begin{align}
U_\ell^{00}&=\left(
\begin{matrix}
2h_\ell\\
2h_\ell\\
(1+(-1)^\ell)h_\ell
\end{matrix}
\right), &
U_\ell^{21}&=\left(
\begin{matrix}
2h_\ell\\
2h_{-\ell}\\
0
\end{matrix}
\right),\\
\intertext{for any $\ell$ and}
U_\ell^{31}&=\left(
\begin{matrix}
2h_\ell\\
0\\
2h_{-\ell}
\end{matrix}
\right), &
U_\ell^{32}&=\left(
\begin{matrix}
0\\
2h_\ell\\
2h_{-\ell}
\end{matrix}
\right)
\label{eq:Ul_ssd}
\end{align}
\end{subequations}
for $\ell$~even and $U_\ell^{31}=U_\ell^{32}=0$ for $\ell$~odd.
Thus by choosing
\begin{equation}
X(\phi)=\sum_{\ell \in \mathbb Z\sm\{0\}} \left(\frac{i}{2\ell}U_\ell^{31}e^{i\ell(\phi_3-\phi_1)}+\frac{i}{2\ell}U_\ell^{32}e^{i\ell(\phi_3-\phi_2)}\right),
\label{X_SSD}
\end{equation}
the dynamics on the perturbed torus
\begin{equation}
f_1^\SSD(\phi)=\sum_{\ell\in\mathbb Z} \left(U_\ell^{00} + U_\ell^{21}e^{i\ell(\phi_2-\phi_1)}\right)
\label{eq:f_all_SSD}
\end{equation} 
is in normal form.
Thus, up to first order, the dynamics is determined by
$f^\SSD_{\delta} = \Omega^\SSD + \delta f_1^\SSD + \mathcal{O}(\delta^2)$ 
and the perturbed torus is embedded through
$e^{\SSD}_{\delta} = E^\SSD+\delta e^{\SSD}_1 + \mathcal{O}(\delta^2)$ 
with $e^{\SSD}_1(\phi) = RX(\phi)+NY(\phi)$, which is the main claim of Proposition~\ref{prop:Comp} for~$\SSD$.

For the biharmonic coupling function
$h(\varphi) = \sin(\varphi+\alpha)+r\sin(2(\varphi+\beta))$
in Proposition~\ref{prop:CompWCoup}, 
we obtain the first-order phase dynamics
\begin{subequations}\label{eq:nf_SSD}
\begin{align}
\dot{\phi}_1&=2\delta(\sin\alpha+r\sin(2\beta)+\sin(\phi_2-\phi_1+\alpha)+r\sin(2(\phi_2-\phi_1+\beta)))\\
\dot{\phi}_2&=2\delta(\sin\alpha+r\sin(2\beta)+\sin(\phi_1-\phi_2+\alpha)+r\sin(2(\phi_1-\phi_2+\beta)))\\
\dot{\phi}_3&=-2+2r\delta\sin(2\beta)
\end{align}
\end{subequations}
Note that in contrast to~$\SDD$, for the choice of this particular perturbation both first and second harmonics contribute to the reduced dynamics on the perturbation of~$\SSD$.
The phase dynamics is embedded in~$\T^6$ via $e^{\SSD}_{\delta} = E^\SSD+\delta e^{\SSD}_1 + \mathcal{O}(\delta^2)$
where
\begin{equation}\label{eq:emb_SSD}
e_1^\SSD(\phi)=\left(
\begin{matrix}
\frac12r(\cos(2(\phi_3-\phi_1+\beta))) \\[0.2cm]
\frac12r(\cos(2(\phi_3-\phi_1+\beta))) \\[0.2cm]
\frac12r(\cos(2(\phi_3-\phi_2+\beta))) \\[0.2cm]
\frac12r(\cos(2(\phi_3-\phi_2+\beta))) \\[0.2cm]
  -\frac12r(\cos(2(\phi_1-\phi_3+\beta))+\cos(2(\phi_2-\phi_3+\beta)))-as(\phi)-bz(\phi) \\[0.2cm]
  -\frac12r(\cos(2(\phi_1-\phi_3+\beta))+\cos(2(\phi_2-\phi_3+\beta)))+as(\phi)+bz(\phi) 
\end{matrix}
\right)
\end{equation}
with
\begin{align*}
a&= \dfrac25(\cos\alpha-2r\sin\alpha), &
s(\phi)&=\cos(\phi_1-\phi_3)+\cos(\phi_2-\phi_3),\\ 
b&=-\dfrac25(\sin\alpha+2r\cos\alpha), &
z(\phi)&=\sin(\phi_1-\phi_3)+\sin(\phi_2-\phi_3).
\label{eq:const_Y_SSD}
\end{align*}
This is the main claim of Proposition~\ref{prop:CompWCoup} for~$\SSD$ with $\varphi = \phi_2-\phi_1$.

%%%%%%
\section{Numerical continuation of solutions on perturbed tori}
\label{sec:dyn_tori}

We now explore numerically how solutions on the perturbed tori change and bifurcate as the asymmetry parameter~$\delta$ is increased.
The approximations computed in the previous section are the key to generate starting data for the numerical continuation: The first-order dynamics allow one to compute approximate solutions for~$\delta$ small on the tori $\SDD, \SSD$.
The embedding then maps these solutions into the full phase space~$\T^6$ to generate starting data that can be continued using~\textsc{Auto} (see~\cite{Code} for accompanying computer code).

Note that the $S_2$~symmetry on~$\T^6$ that acts by permuting the oscillators within one population induces a symmetry action on~$\T^3$. If a population is desynchronized~(`$\Dp$') this corresponds to a shift by~$\pi$.
We expect that this is reflected in the bifurcation diagrams.
Moreover, this explains why a perturbation with single harmonics does not affect the local dynamics on~$\SDD$ where only even harmonics appear in the normal form dynamics: 
Only even harmonics preserve a phase shift by~$\pi$.

\begin{figure}%[htp]
\begin{center}
  % replace aims_logo.pdf by your figure file name
  \includegraphics[width=\linewidth]{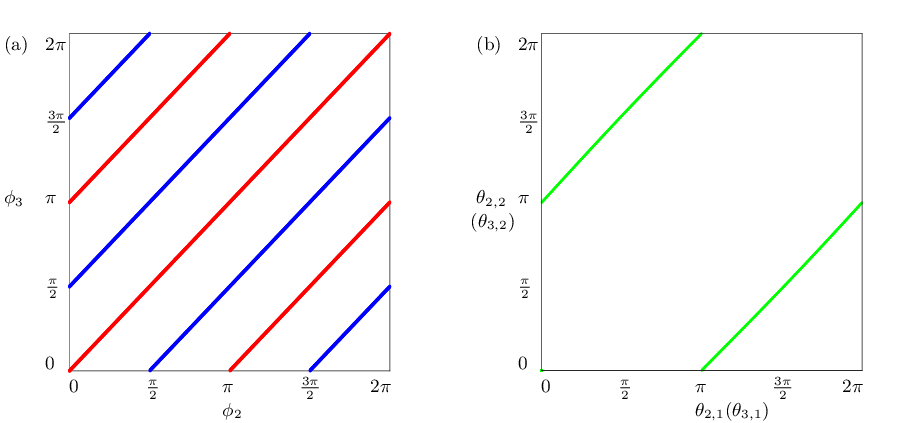}
  \caption{
Dynamics of system~\eqref{eq:nf_SDD} for $\alpha=\beta=\frac{\pi}{2}$, $r=0.2$ and $\delta=0.01$. 
Panel~(a) shows solutions to the reduced system~\eqref{eq:SDD_red} on the torus~$\T^2$ described by coordinates $(\phi_2,\phi_3)$ representing the `$\Dp$' populations; these periodic orbits wind around~$\T^2$. Shown are
two stable periodic orbits (blue) at $\psi\in\{\frac{\pi}{2}, \frac{3\pi}{2}\}$ and two unstable periodic orbits (red) at $\psi\in\{0, \pi\}$. Panel~(b) shows a periodic orbit in the full system that lies on $\SDD$ for $\phi_3-\phi_2 = 0$. Note that the two `$\Dp$' populations have approximately the same phase.}\label{fig:nf_sdd}
  \end{center}
\end{figure}

%%%
\subsection{Dynamics on the perturbation of $\SDD$}
\label{sec:dyn_sdd}

To understand the dynamics of system~\eqref{eq:Dynamics3x2} on~$\SDD$ with the perturbation driven by~\eqref{eq:Zj} and~\eqref{eq:h2} as the symmetry breaking parameter~$\delta$ is varied, we first consider the first-order approximation of the dynamics in~$\delta$.
With $\psi=\phi_3-\phi_2$ being the phase difference between the two desynchronized populations, the dynamics up to first order are given by~\eqref{eq:SDD_red_lem} which read
\begin{equation}
\begin{split}
\dot \phi_1 &= 2\delta(\sin\alpha+r\sin(2\beta)),\\
\dot \psi &= 4r\delta(\sin (2(-\psi+\beta))-\sin(2(\psi+\beta))).
\end{split}
\label{eq:SDD_red}
\end{equation}
As~$\phi_1$ has constant motion, the dynamics are effectively one-dimensional; the $\psi$-nullclines correspond to periodic orbits in the full system~\eqref{eq:nf_SDD}, with their stability given by the equation for $\dot\psi$. 
For $\delta,r>0$, the dynamics depends on~$\beta$ as follows: 
For $\beta^*=\frac{\pi}{4}+\frac{k\pi}{2}, k\in\Z,$ the torus~$\T^2$ is foliated by neutrally stable periodic orbits. 
For $\beta\neq \beta^*$, the system possesses four hyperbolic periodic orbits, two of which are stable and two are unstable; this is a consequence of the second harmonics in the normal form. 
The two stable (resp.~unstable) periodic orbits are related by a shift by~$\pi$ along the corresponding population (i.e., $\phi$-coordinate) and have exactly the same Floquet multipliers (computed numerically), as expected from the $\mathbb Z_2$-symmetry translated into normal form coordinates. 
The stability of these periodic orbits is reversed once~$\beta$ passes through~$\beta^*$.
The location of these hyperbolic periodic orbits does not change upon variations of $\beta\neq\beta^*$. 
We fix $\beta=\frac{\pi}{2}$ so that~\eqref{eq:SDD_red} has unstable periodic orbits for $\psi\in\{0,\pi\}$ and {stable} periodic orbits for $\psi\in\big\{\frac{\pi}{2}, \frac{3\pi}{2}\big\}$.
They correspond to the periodic orbits for which $\phi_3-\phi_2\in\{0,\pi\}$ and $\phi_3-\phi_2\in\big\{\frac{\pi}{2},\frac{3\pi}{2}\big\}$, respectively; this gives an idea on how the `$\Dp$' populations organize relative to each other. 
\Fref{fig:nf_sdd} shows the organization of phase space for system~\eqref{eq:SDD_red} with $\beta=\frac{\pi}{2}$. Panel~(a) shows the configuration in $(\phi_2,\phi_3)$-plane, relating the `$\Dp$'~populations.
Panel~(b) shows the unstable periodic orbit that corresponds to $\phi_3-\phi_2=0$ in the full system.

\newcommand{\Auto}{{\sc Auto}}

We now want to study the dynamics on~$\SDD$ in full space~$\T^6$ for different values of~$\delta$ using numerical continuation. 
To generate initial data, we start with an asymptotically stable periodic obit of the first-order approximation~$f_1^\SSD$ on~$\mathbb{T}^3$.  
These periodic orbits can then be continued by taking advantage of the collocation framework of \Auto{}.
Once we obtain a good approximation of one of the stable periodic orbits corresponding to $\psi=\frac{\pi}{2}$ or $\psi=\frac{3\pi}{2}$ in system~\eqref{eq:nf_SDD}, we use the embedding~\eqref{eq:emb_SDD} to get a periodic orbit of the perturbed system~\eqref{eq:Dynamics3x2} on~$\T^6$ that lies on a perturbation of~$\SDD$. 
Note that such a periodic orbit is not stable in~$\T^6$ since we carry on instabilities of the saddle torus~$\SDD$, whose stability properties persist under small perturbations due to normal hyperbolicity.
This is one of the benefits of  the reduction to the perturbed tori: 
We directly obtain approximations for two (symmetry-related) saddle periodic orbits, which is practically impossible by standard shooting techniques.
We use any of the periodic orbits obtained this way as an initial condition for studying numerically the dynamics on a perturbation of~$\SDD$ with respect to~$\delta$. 

\begin{figure}%[htp]
\begin{center}
  % replace aims_logo.pdf by your figure file name
  \includegraphics[width=\linewidth]{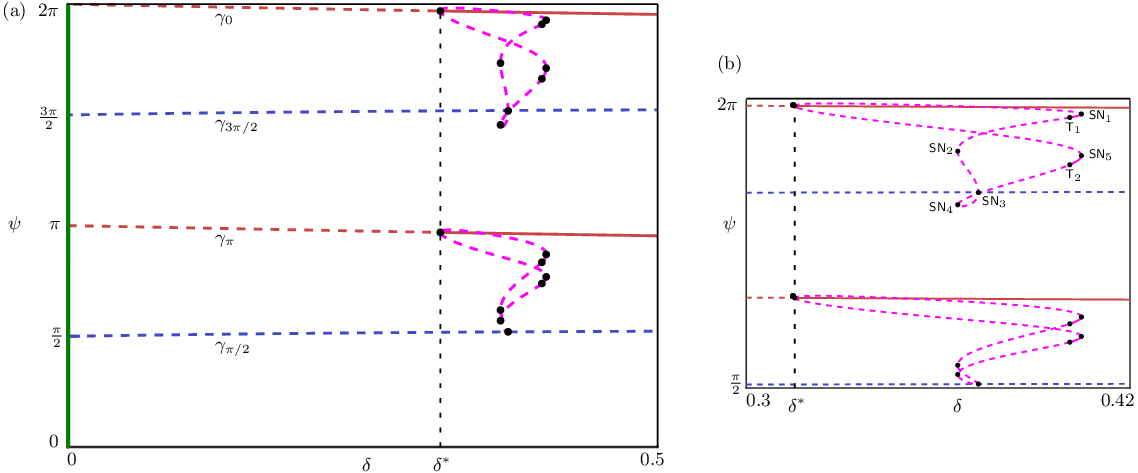}
  \caption{
Bifurcations can stabilize periodic orbits of the full system~\eqref{eq:Dynamics3x2} stemming from~$\SDD$ for a perturbation determined by~\eqref{eq:Zj} and~\eqref{eq:h2} that is sufficiently large $\delta>0$. 
The rest of the parameters are as in~\eqref{eq:par}. 
Panel~(a) shows the four branches of hyperbolic periodic orbits $\gamma_0, \gamma_{\frac{\pi}{2}}, \gamma_{\pi}$ and $\gamma_{\frac{3\pi}{2}}$ emanating from $\psi^* \in \big\{0, \frac{\pi}{2},\pi,\frac{3\pi}{2}\big\}$ at $\delta=0$, bifurcations that change their stability in the full system, and secondary solution branches.
The vertical green line at $\delta=0$ represents the unperturbed case, in which the torus~$\SDD$ is foliated by neutrally stable periodic orbits. 
Panel~(b) highlights bifurcations at $\delta=\delta^*$ that stabilize the branches~$\gamma_0$ and~$\gamma_{\pi}$ in the full system and corresponding  secondary branches of periodic orbits forming isolas.}\label{fig:pobranches_sdd}
  \end{center}
\end{figure}

\Fref{fig:pobranches_sdd} shows a bifurcation diagram for the periodic orbits lying on a perturbation of~$\SDD$, for small values of $\delta>0$ in the full system~\eqref{eq:Dynamics3x2}. 
These orbits were obtained from numerical continuation as a boundary-value-problem with \Auto{} using the initial data as described. 
Here we chose the fixed parameter values 
%\JPM{Check parameter $\alpha_2$.}
\begin{subequations}\label{eq:par}
\begin{align}
\alpha_2 &= \frac{\pi}{2}, & \alpha_4&=\pi, &  K&=0.4, & r_0&=0.1
\label{eq:par_unp}
\intertext{that specify the unperturbed dynamics
and}
\alpha&=\frac{\pi}{2}, & \beta &=\frac{\pi}{2}, & r&=0.2
\label{eq:par_pert}
\end{align}
\end{subequations}
that specify the symmetry breaking perturbation determined by~\eqref{eq:h2}.
Parameters~\eqref{eq:par_unp} are such that the torus~$\SDD$ and~$\SSD$ are of saddle type for~$\delta=0$; see~\cite{Bick2019}. 
We set $\delta=0.01$ at the beginning of the computation to obtain the initial stable periodic orbit in system~\eqref{eq:nf_SDD} and then embed it in~$\T^6$ to start the continuation. 
The horizontal axis in \fref{fig:pobranches_sdd} corresponds to the values of the perturbation parameter~$\delta$, and the vertical axis corresponds to the average of the phase difference $\theta_{3,1}-\theta_{2,1}$, modulo $2\pi$, along the corresponding periodic orbit for each parameter value. 
We abuse notation and denote this  averaged phase difference by~$\psi$ as in \fref{fig:nf_sdd}, in order to be able to compare both figures. 
In fact,~$\psi$  measures the relative positions between the two `$\Dp$' populations in~$\SDD$. 
The vertical green line in Panel~(a) represents the family of neutrally stable periodic orbits that exist for $\delta=0$; they correspond to the foliation of~$\SDD$ by periodic orbits in the fully symmetric case.

\begin{figure}%[htp]
\begin{center}
  % replace aims_logo.pdf by your figure file name
  \includegraphics[width=\linewidth]{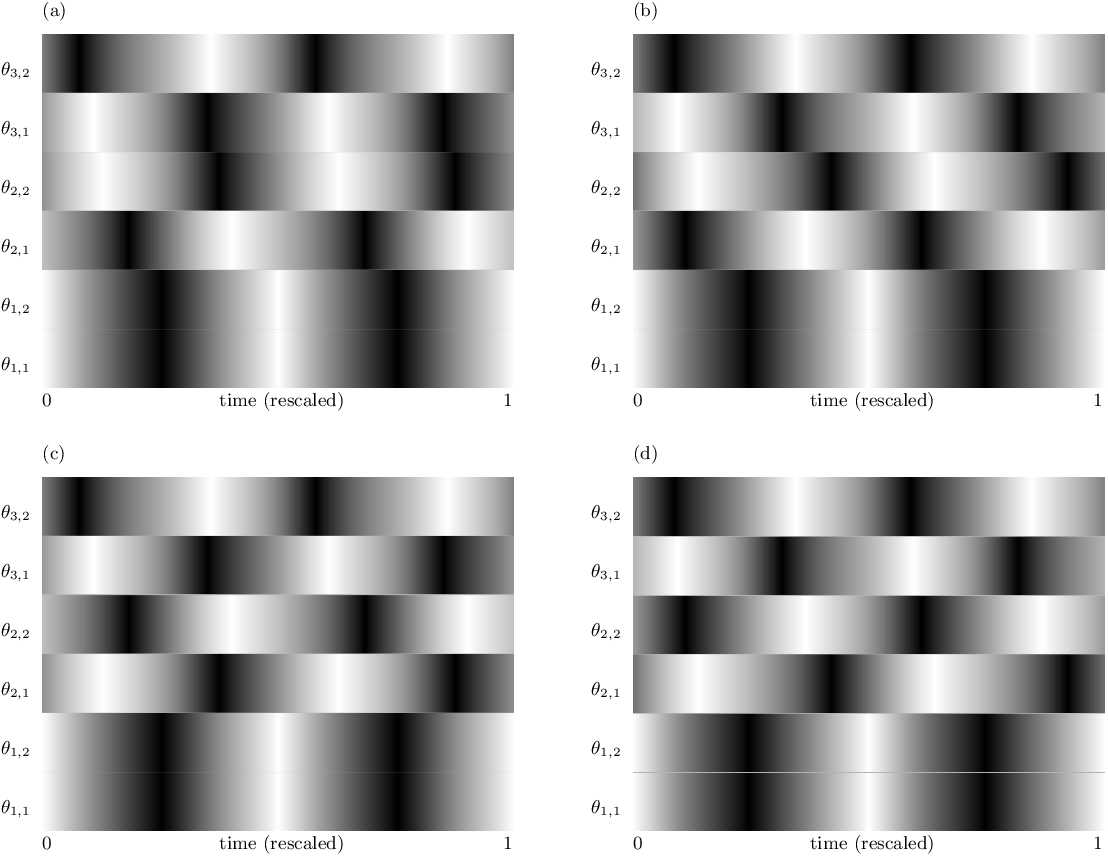}
  \caption{
Stable synchrony patterns bifurcating from primary branches of periodic orbits. 
Panels~(a)--(d) show the phase evolution of oscillators in `$\Dp$' populations along one period on coexisting stable periodic orbits lying on the secondary branches emanating from $\psi=0$ (bottom) and $\psi=\pi$ (top) for $\delta\approx 0.40372$; see the magnified parameter range in \fref{fig:pobranches_sdd}(b). 
Shading indicates the deviation from~$\pi$ of each phase oscillator, where black indicates that $\theta_{\sigma,k}=\pi$ and white $\theta_{\sigma,k}=0$ or~$2\pi$.}\label{fig:st_patterns}
  \end{center}
\end{figure}

As the asymmetry parameter~$\delta$ is varied, four branches of periodic orbits emanate for $\delta>0$ from the vertical green line in \fref{fig:pobranches_sdd}(a) at $\psi^*\in\big\{0, \frac{\pi}{2}, \pi, \frac{3\pi}{2}\big\}$ denoted by $\gamma_{\psi^*}$; see the resemblance with \fref{fig:nf_sdd}, which corresponds to the (reduced) dynamics on phase space for a fixed value of $\delta\neq 0$ small. 
These four periodic orbits can be interpreted in terms of their phase configurations as invariant subsets of~$\SDD$: 
In~$\gamma_0$ the two `$\Dp$' populations are phase synchronized, in~$\gamma_{\frac{\pi}{2}}$ they are $\frac{\pi}{2}$ apart etc---of course this interpretation is only approximate as~$\delta$ is increased.
We use similar colors as for the branches of periodic orbits in \fref{fig:nf_sdd}, again to stress the connection between the two figures, and highlighting the stability within~$\SDD$ for~$\delta$ small (they are unstable in~$\T^6$).
For each~$\delta>0$, the periodic orbits on $\gamma_{0}$ and $\gamma_{\pi}$ have (numerically) the same Floquet multipliers; the same occurs with the periodic orbits on~$\gamma_{\frac{\pi}{2}}$ and~$\gamma_{\frac{3\pi}{2}}$. 
This is reminiscent of the~$\Z_2$ symmetry in the (unperturbed) full system that is still present in the full system for $\delta\neq 0$ small.

Note that the branches~$\gamma_0$ and~$\gamma_{\pi}$ both become stable at (numerically) the same value~$\delta=\delta^*\approx 0.3152$, from which a secondary branch of periodic orbits emanates on each primary branch $\gamma_{0}, \gamma_{\pi}$.
The secondary bifurcations are shown in detail in Panel~(b) of \fref{fig:pobranches_sdd}.
Each of these secondary branches forms an isola of periodic orbits, and can be thought of as two secondary branches emanating from the main branch at~$\delta^*$, that collide at the point~$\mathsf{SN}_3$. 
The two secondary branches forming each isola undergo exactly the same bifurcations at (numerically) exactly the same parameter values: both (secondary) branches pass through a saddle-node~($\mathsf{SN}_{1,5}$), torus~($\mathsf{T}_{1,2}$) and another saddle-node bifurcation~($\mathsf{SN}_{2,4}$) before colliding on~$\mathsf{SN}_3$. 
This suggests that the secondary branches arise from a pitchfork{-like} bifurcation whose symmetries are not present in the projection used in \fref{fig:pobranches_sdd}. 
This bifurcation must be transverse to the perturbed torus~$\SDD$ (if it still exists), since all the Floquet multipliers of the corresponding periodic orbits are stable after the bifurcation.
This suggests that the pitchfork-like scenario is combined with a different phenomenom. 
In fact, the secondary branches arising from this bifurcation cross the branch~$\gamma_{\frac{3\pi}{2}}$, indicating that a torus breakdown has likely happened for smaller values of~$\delta$. 
This is consistent with the existence of saddle-node bifurcations on the secondary branches before this crossing, which often relates to a loss of smoothness of the torus. 
A deeper exploration of the bifurcations involving the torus breakdown is beyond of the scope of this paper, as we focus on what emanates from this primary bifurcation scenario.
There is a small interval of multistability between the points~$\mathsf{SN}_1$ ($\mathsf{SN}_5$) and~$\mathsf{T}_1$ ($\mathsf{T}_2$), where two stable periodic orbits on each secondary branch coexists with a stable periodic orbit lying on the primary branch, adding up to 6 stable periodic orbits all together. 
A similar structure occurs along the primary branch $\gamma_{\pi}$, but the labels are not included in \fref{fig:pobranches_sdd}(b).

The small interval of multistability due to the existence of bifurcating secondary branches leads to new coexisting stable synchrony patterns. 
\Fref{fig:st_patterns} shows the phases of the oscillators along one period of these secondary stable periodic orbits. 
Shown are the phases of the six oscillators against the rescaled integration time for $\delta \approx 0.403722$, that is within the interval between a saddle-node and a torus bifurcation where the periodic orbits on the isolas on \fref{fig:pobranches_sdd} are all stable. The color indicates the deviation from~$\pi$ of~$\theta_{\sigma,k}$, where black means $\theta_{\sigma,k}=\pi$ and white means $\theta_{\sigma,k}=0$ or $2\pi$. 
Here, Panels~(a) and~(b) correspond to stable periodic orbits lying on top and bottom of the secondary branch that emanates from $\psi=\pi$, respectively and Panels~(c) and~(d) show the same for the the secondary branch emanating from~$\psi=0$. 
This shows that in these new synchrony patterns  the phases on `$\Dp$' populations are not exactly~$\pi$ apart, but rather wiggle around~$\pi$ and even get close for extremely short periods of time.

\begin{figure}%[htp]
\begin{center}
  % replace aims_logo.pdf by your figure file name
  \includegraphics[width=\linewidth]{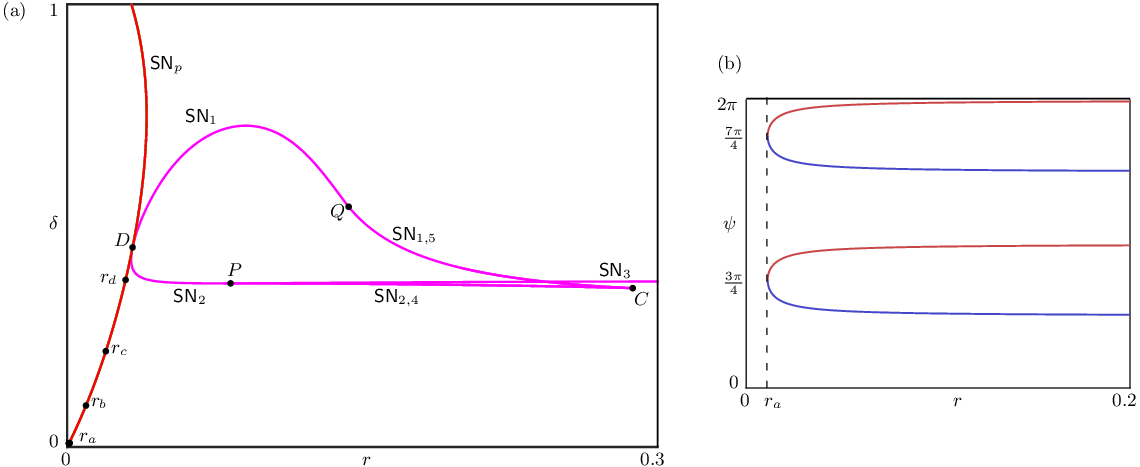}
  \caption{
Some codimension-two bifurcations of system~\eqref{eq:Dynamics3x2} with the perturbation driven by~\eqref{eq:Zj} and~\eqref{eq:h2}. 
Panel~(a) shows curves of saddle-node bifurcation arising on the primary (red) and secondary (magenta) branches that organize the parameter plane~$(r,\delta)$. 
See the main text for an elaboration on the  codimension-two points shown.
Panel~(b) shows the origin of the primary saddle-node bifurcation as one continues a primary branch of periodic orbits in parameter~$r$, with $r$~decreasing. 
This bifurcation gives rise to the curve~$\mathsf{SN}_p$ in Panel~(a).}\label{fig:R_delta_sdd}
  \end{center}
\end{figure}

The organization of the saddle-node bifurcations along the secondary branches shown in \fref{fig:pobranches_sdd} can be understood in terms of a two-parameter continuation in both~$\delta$ and the strength~$r$ of the second harmonic of the symmetry breaking perturbation~\eqref{eq:h2}.
\Fref{fig:R_delta_sdd}(a) shows curves of saddle-node bifurcations (magenta) that correspond to $\mathsf{SN}_1$--$\mathsf{SN}_5$ on the secondary branch of periodic orbits emanating from~$\gamma_0$, on the parameter plane $(r,\delta)$. 
The isola shown in \fref{fig:pobranches_sdd}(b) is for $r=0.2$, where there are saddle-node bifurcations at three different~$\delta$ values and with~$\mathsf{SN}_1$ and~$\mathsf{SN_5}$ (and also $\mathsf{SN}_2$ and $\mathsf{SN_4}$) occurring simultaneously. 
When continued in~$r$ and~$\delta$, the curves~$\mathsf{SN}_1$ and~$\mathsf{SN}_5$ are on top of each other until one of them disappears. 
The same occurs with the curves~$\mathsf{SN}_2$ and $\mathsf{SN}_4$.
Summarizing the results, we observe different codimension-two points:
\begin{itemize}
\item The point~$C$ corresponds to a (double) cusp bifurcation, in which the curve~$\mathsf{SN}_1$ ($\mathsf{SN}_5$) and~$\mathsf{SN}_2$ ($\mathsf{SN}_4$) collide and disappear, while~$\mathsf{SN}_3$ persists. 
In this process, there is a point in~$(r,\delta)$ for which~$\mathsf{SN}_1$, $\mathsf{SN}_5$ and $\mathsf{SN}_3$ occur simultaneously and then swap places with~$\mathsf{SN}_3$.
\item At the point~$Q$, the curve~$\mathsf{SN}_5$ terminates and~$\mathsf{SN}_1$ continues.
\item Similarly, at the point~$P$, the curve~$\mathsf{SN}_4$ terminates and~$\mathsf{SN}_2$ continues. 
This is the point in parameter plane from which the curve~$\mathsf{SN}_3$ emerges.
\item The point~$D$ corresponds to a degenerate codimension-two point, in which~$\mathsf{SN}_1$ and~$\mathsf{SN}_2$ collide and terminate. 
At the point~$D$ these two curves come tangent to the saddle-node curve~$\mathsf{SN}_p$ (red curve).
\end{itemize}
A similar structure is expected for the saddle-node bifurcations occurring on the secondary branch of periodic orbits emanating from~$\gamma_{\pi}$.

\begin{figure}%[htp]
\begin{center}
  % replace aims_logo.pdf by your figure file name
  \includegraphics[width=\linewidth]{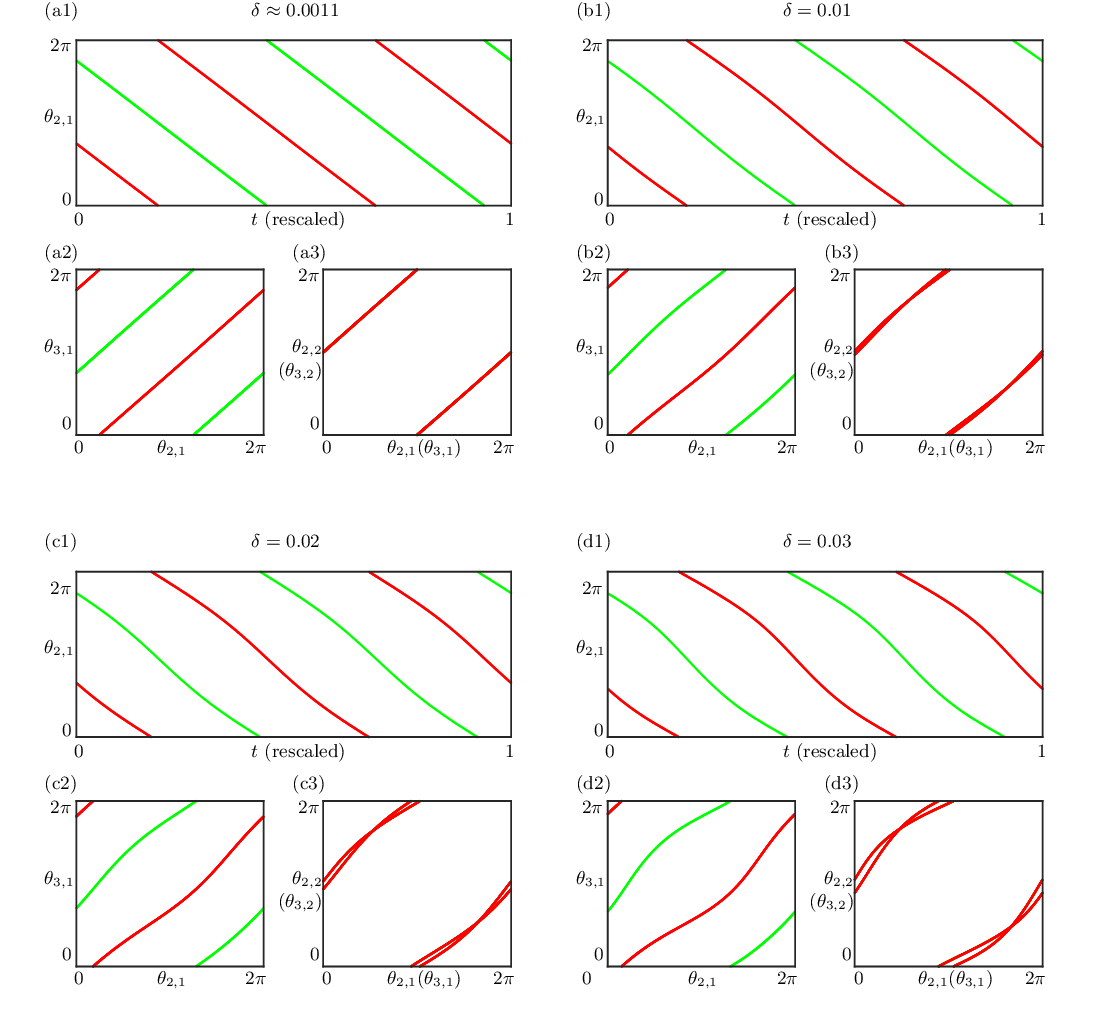}
  \caption{
Evolution of periodic orbits along the overlapping principal saddle-node curves~$\mathsf{SN}_p$ in \fref{fig:R_delta_sdd}.
Panels (a)--(d) correspond to the points $r_a,r_b, r_c, r_d$, respectively.
Each panel shows (1)~a time series of the coordinate~$\theta_{2,1}$ with respect to the rescaled \Auto ~integration time for orbits on each principal saddle-node curve in green and red, (2)~a relative position between the two `$\Dp$' populations on each principal saddle-node curve in green and red, and (3)~how each `$\Dp$'  populations deviate from~$\pi$ as $(r,\delta)$ moves away from the origin.}\label{fig:evo_sn_sdd}
  \end{center}
\end{figure}

The saddle-node curve $\mathsf{SN}_p$ in \fref{fig:R_delta_sdd}(a) corresponds to a saddle node bifurcation of the primary branches of periodic orbits~$\gamma_{\frac{\pi}{2}}$ (dark blue) and~$\gamma_{\pi}$ (dark red) as~$r$ is decreased at a value of~$r$ that is~$\mathcal O(\delta^{3/2})$; cf.~\Fref{fig:R_delta_sdd}(b).
(The same occurs with the branches~$\gamma_{\frac{3\pi}{2}}$ (dark blue) and $\gamma_{0}$ (dark red) at numerically the same value.)
Starting with $r=0.2$ and decreasing~$r$ while keeping $\delta=0.01$ fixed, the saddle-node bifurcation occurs at $r=r_a\approx 0.0011$. 
At these bifurcation points, the phases of the `$\Dp$'~populations of are~$\frac{3\pi}{4}$ and~$\frac{7\pi}{4}$ apart, which is half way between the relative positions of the $`\Dp$'~populations initially. 
This saddle-node bifurcation is a consequence of the~$\mathbb Z_2$ symmetry inherited from reduced dynamics, since by decreasing~$r$ towards zero we are decreasing the second harmonic in the symmetry breaking perturbation that allow the solutions to exist in the first place.
This saddle-node bifurcation on the primary branches can be then continued in~$(r,\delta)$ tracing out the curve~$\mathsf{SN}_p$. 
The curve converges to~$(0,0)$ in parameter plane; solutions along~$\mathsf{SN}_p$ at $r_a\approx 0.0011$, $r_b = 0.01$ (panel (b)), $r_c = 0.02$, and~$r_d = 0.03$ are shown in \Fref{fig:evo_sn_sdd}.
Naturally, the curve~$\mathsf{SN}_p$ bounds the existence of other bifurcations $\mathsf{SN}_1$--$\mathsf{SN}_5$ along the secondary branches since the primary branches $\gamma_{\psi^*}$ with $\psi^*\in\{0,\frac{\pi}{2}, \pi, \frac{3\pi}{2}\}$ emerge there.

%%%
\subsection{Dynamics on the perturbation of~$\SSD$}
\label{sec:dyn_ssd}

The first-order approximation of the dynamics on~$\SSD$ are given by~\eqref{eq:SSD_red_lem}, which read
\begin{equation}
\begin{split}
\dot \varphi &= -4\delta\sin\varphi(\cos \alpha+2r\cos\varphi\cos(2\beta)),\\
\dot \phi_3 &= -2+2r\delta\sin(2\beta)
\end{split}
\label{eq:SSD_red}
\end{equation}
with $\varphi=\phi_2-\phi_1$.
The variable~$\phi_3$ has constant motion and the nullclines of~$\varphi$ provide periodic orbits in system~\eqref{eq:nf_SSD}.
The quantity $\rho=\rho(r,\alpha,\beta)=\frac{\cos\alpha}{2r\cos(2\beta)}$ determines the influence of second harmonics in the $\varphi$-dynamics unless $r=0$ or $\beta=\beta^*$ with $\cos(2\beta^*)=0$.
For $\lvert \rho \rvert>1$ there is one stable and one unstable periodic orbit located at $\phi=0,\pi$; which of these orbits is stable/unstable depends on the sign of~$\rho$. 
The periodic orbit $\varphi=\pi$ undergoes a pitchfork bifurcation at $\lvert \rho \rvert =1$ and two symmetry-related periodic orbits~$\varphi=\varphi^*$ and~$\varphi=2\pi-\varphi^*$ emanate from there and persist for $\lvert \rho \rvert <1$, where $\varphi^*$ is solution of $\cos\alpha+2r\cos\varphi\cos(2\beta)=0$. 
The critical value $\lvert \rho \rvert =1$ is a threshold for the onset of the influence of second harmonics in the dynamics of \eqref{eq:nf_SSD}. 
For $\lvert \rho \rvert <1$ system~\eqref{eq:nf_SSD} has two stable and two unstable periodic orbits and the phase space looks similar to \fref{fig:nf_sdd}, which is the situation for the parameter values given by~\eqref{eq:par}. 
Note that the choice of~$\varphi$ in the reduced system~\eqref{eq:SSD_red} has a different interpretation than the variable~$\psi$ for~$\SDD$: 
Rather than the relative position between desynchronized populations, $\varphi$~corresponds the the relative position between the synchronized (`$\Sp$') populations on~$\SSD$. 

\begin{figure}%[htp]
\begin{center}
  % replace aims_logo.pdf by your figure file name
  \includegraphics[width=\linewidth]{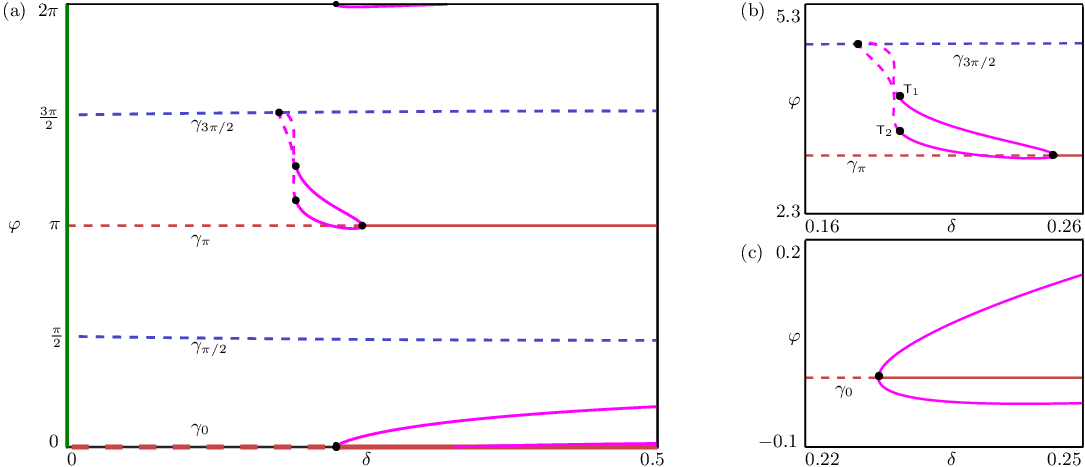}
  \caption{
Bifurcations of periodic orbits of the full system~\eqref{eq:Dynamics3x2} with the perturbation driven by~\eqref{eq:Zj} and~\eqref{eq:h2} for $\delta>0$ for the perturbed torus~$\SSD$. 
The rest of the parameters are as in \fref{fig:pobranches_sdd}. 
Panel~(a) shows four branches of hyperbolic periodic orbits $\gamma_0, \gamma_{\frac{\pi}{2}}, \gamma_{\pi}, \gamma_{\frac{3\pi}{2}}$ emanating from $\psi^*\in \big\{0, \frac{\pi}{2},\pi, \frac{3\pi}{2}\big\}$ at $\delta=0$. 
The vertical green line located at $\delta=0$ represents the unperturbed case, in which the torus $\SSD$ is foliated by neutrally stable periodic orbits. 
Panel~(b) highlights the isola of secondary branches of periodic orbits emanating from the branch~$\gamma_{\pi}$, and Panel~(c) shows a zoom of the origin of the branch of secondary periodic orbits that emanates from~$\gamma_{0}$.}\label{fig:pobranches_ssd}
  \end{center}
\end{figure}

We use the first-order approximation in the same way as in \sref{sec:dyn_sdd} to generate starting data for numerical continuation. 
Parameter values here are fixed and given by~\eqref{eq:par}, so that the periodic orbit corresponding to $\varphi\in\{\frac{\pi}{2},\frac{3\pi}{2}\}$ are stable and $\varphi\in\{0,\pi\}$ are unstable.
This way we numerically obtain an initial periodic orbit on~$\T^6$ that can be continued in~$\delta$.
For $\delta=0$ (the unperturbed system) the torus~$\SSD$ is foliated by neutrally stable periodic orbits.
Four branches of saddle periodic orbits bifurcate from this family at $\varphi\in\{0,\frac{\pi}{2}, \pi, \frac{3\pi}{2}\}$, and persist for $\delta\neq0$ small; we denote this branches by $\gamma_0,\gamma_{\frac{\pi}{2}},\gamma_{\pi}$ and $\gamma_{\frac{3\pi}{2}}$ in analogy to those on~$\SDD$ in~\sref{sec:dyn_sdd}.

\Fref{fig:pobranches_ssd} shows the bifurcations of these periodic orbits for small values of $\delta>0$.
The branch~$\gamma_0$ becomes stable at $\delta\approx0.227$, where two secondary branches of unstable periodic orbits (magenta) emerge from the bifurcation point; cf.~\Fref{fig:pobranches_ssd}(c).
In contrast to \fref{fig:pobranches_sdd}, these two secondary branches do not connect back forming an isola, but rather remain unstable.
The branch~$\gamma_{\frac{\pi}{2}}$ remains unstable and undergoes a saddle-node bifurcation that connects it with a branch that emanates from a homoclinic orbit at $\delta\approx 0.796$ (not shown). 
The branch~$\gamma_{\pi}$ gains stability at $\delta=0.25$, and then terminate at $\delta\approx 1.67$. 
Two secondary branches of stable periodic orbits (magenta) emanate from the bifurcation point and coexist with the primary branch for $\delta<0.25$. These two secondary branches undergo the same bifurcations at (numerically) the same parameter values as for~$\SDD$ (cf.~\fref{fig:pobranches_sdd}), but with the difference that the branches here connect with the primary branch~$\gamma_{\frac{3\pi}{2}}$ in a branching point on~$\gamma_{\frac{3\pi}{2}}$ at $\delta\approx 0.179$. 
Starting from $\delta=0.25$ and decreasing~$\delta$, these secondary branches undergo a torus bifurcation to become unstable, then pass through a saddle-node bifurcation to finally meet at the primary branch~$\gamma_{\frac{3\pi}{2}}$; cf.~\Fref{fig:pobranches_ssd}(b). 
Finally, the branch $\gamma_{\frac{3\pi}{2}}$ has two more branching  points before terminating at $\delta\approx 0.85$, where the periodic orbits approach heteroclinic trajectories.

%%%%%%
\section{Discussion}
\label{sec:discussion}

We analyzed the dynamics on invariant tori that arise by perturbing relative equilibria in the context of a network of coupled oscillators. 
First, we showed that the parameterization method developed in~\cite{vondergracht2023parametrisation} can be applied to relative equilibria of a continuous~$\T^m$ symmetry.
This yields explicit equations for the perturbed tori as well as the dynamics thereon to any order.
Second, we computed the first-order approximation of the torus and its dynamics explicitly by solving the corresponding conjugacy equations.
Third, we used the first-order approximation to generate starting data for numerical continuation in \Auto. 
This allowed to compute branches of periodic orbits that emerge on the perturbed tori as the coupling strength is increased beyond the regime where the first-order approximation is valid (or even beyond the existence of an invariant torus).
The choice of explicitly parameterized phase interaction functions yields a link to parameters in physical models, e.g., through phase reduction~\cite{Bick2023,Nakao2015,Leon2019a}.

Our findings provide insight into the effect of forced symmetry breaking on the local dynamics of invariant tori existing in the coupled oscillator network~\eqref{eq:Dynamics3x2}. 
There are symmetries on these sets that relate to synchrony patterns and make the local dynamics on the tori unaffected by the contribution of particular higher harmonics in a perturbation.
Importantly, for $\delta>0$ large enough we stabilize new synchrony patterns via secondary bifurcations that are a result of forced symmetry breaking in the system.
Indeed, for sufficiently large parameter values we observe bifurcations leading to multistability that involves periodic orbits that do not exist without broken symmetry.

Note that the numerical continuation focuses on periodic solutions, i.e., synchrony patterns, on an invariant torus that persists for small~$\delta$ but not on the torus itself. Hence, the periodic solutions may continue to exist beyond a torus breakdown. Indeed, we find numerical evidence for a torus breakdown to happen before a branch of periodic solutions ceases to exist; cf.~\sref{sec:dyn_sdd}.

Understanding global collective dynamics of phase oscillator networks remains an exciting challenge. For some fixed values of the parameters~$\alpha_2$ and~$\alpha_4$, the unperturbed  network described by~\eqref{eq:Dynamics3x2} with coupling functions~\eqref{eq:coupling} supports robust heteroclinic dynamics of localized frequency synchrony~\cite{Bick2019}.
Here the invariant tori~$\SDD$,~$\SSD$ of saddle type and their images under the~$\mathbb Z_3$-action form a heteroclinic structure
\begin{equation*}
\mathrm{DSS} \leadsto \mathrm{DDS} \leadsto \mathrm{SDS} \leadsto \mathrm{SDD} \leadsto \mathrm{SSD} \leadsto \mathrm{DSD} \leadsto \mathrm{DSS},
\label{eq:hcycle}
\end{equation*}
which  exists for an open set of parameters~$(q,K)$. 
The work done in this paper can be used in the larger problem of studying global features of system~\eqref{eq:Dynamics3x2} that relate to perturbations of the heteroclinic structure. 
We can now use the knowledge gathered about the local properties of the tori $\SSD$, $\SDD$, etc.~and use it to study how the connections between these sets are affected under forced symmetry breaking:
One would expect that trajectories close to the stable manifold of a perturbed saddle torus~$\Tbf$ first approach the torus, then move towards a saddle periodic orbit~$\gamma$ on~$\Tbf$ (at a timescale determined by the perturbation parameter~$\delta$), before leaving the neighborhood of~$\Tbf$ near the unstable manifold of~$\gamma$.
Thus, the unstable manifold of~$\gamma$ (within the unstable manifold of~$\Tbf$) gives information about the global dynamics of the perturbed system.
While this is a question for further research, the numerical setup implemented here---using continuation of orbit segments in a boundary-value-problem setup---can be also used, together with state-of-the-art techniques to compute invariant manifolds of saddle invariant objects, to tackle some of these global questions. For an entry point on these continuation methods we refer the reader to \cite{Redbook}.

\section*{Code Availability}
\noindent
\textsc{Auto} code implementing the numerical continuation is available on GitHub~\cite{Code}.

%%%%%%%%%%%%%%%%%%%%%%%%%%%%%%%%%%%%%%%%%%%%%%%%%%%%%%
%          7. REFERENCES SECTION
%%%%%%%%%%%%%%%%%%%%%%%%%%%%%%%%%%%%%%%%%%%%%%%%%%%%%%

%       READ THIS SECTION CAREFULLY

% Each of the references below MUST be cited in your article above. Do not include references that are not cited in your article.

% Follow the examples below carefully. We strongly suggest that you copy and paste your reference information directly into our examples.

% List all references in alphabetical order according to the first author's last name.

% Verify each URL works correctly and can be accessed properly. Your URL links should be to reputable websites. The command line for a website link begins with: \url{ }

% Do not add MR or DOI numbers to your references. AIMS production staff will add this information.

% Using BibTex is not recommended but can be handled.

\bibliographystyle{unsrt}
\bibliography{bibfile}

% APPENDIX

%%%%%%
\section*{Appendix: Numerical implementation}
Here we briefly discuss the numerical setup used for the computation of the periodic orbits and bifurcation diagrams of \sref{sec:dyn_tori}. This is done via continuation of solutions to a two-point boundary value problem implemented in the package {\sc Auto}~\cite{auto}; see~\cite{Code} for \textsc{Auto} code for this paper.  
In contrast to standard shooting methods, the continuation routines of {\sc Auto} use orthogonal collocation with piecewise polynomials~\cite{Ascher,Boor}, and the size of the pseudo-arclength continuation step is determined from the entire orbit segment. This computational approach copes very well with sensitive systems and with systems defined on a torus; see~\cite{Redbook} for more background information.

As standard in {\sc Auto}, we rescale time and write the system  in the form
\begin{equation}
u^{\prime}=TF(u,\lambda).
\label{eq:rescaled}
\end{equation}  
Here, $u=(u_1,\dotsc, u_6) \in \T^6$, $F=F_{\delta}:\T^6 \times \mathbb R^k \to \mathbb \T^6$ is the right-hand side of~\eqref{eq:Dynamics3x2} and $\lambda \in \mathbb R^k$ is its vector of parameters, including the perturbation parameter~$\delta$. Importantly, $u:[0,1]\to\T^6$ so that any orbit segment is parameterized over the unit interval~$[0,1]$; the actual integration time~$T$ is considered as a separate parameter. The function~$u$ is a unique solution of~\eqref{eq:rescaled} if suitable boundary conditions are imposed at one or both end points~$u(0)$ and~$u(1)$. Therefore, each orbit segment is defined in terms of the conditions one imposes upon~$u(0)$ and~$u(1)$.

The usual boundary condition for a periodic in $\mathbb R^n$ is $u(1)-u(0)=0$. This is already implemented in the {\sc Auto} routines for the continuation of a periodic orbit in parameters when we set the corresponding problem type by choosing {\sc IPS}=2 in the {\sc Auto}-constants file; see~\cite{auto}. In the context of coupled oscillators, since the phase variable~$\theta$ (and therefore its rescaled version~$u$) lies on a torus a periodic orbit can be defined by
\begin{equation}
u_i(1)-u_i(0)-2k_i\pi=0, \quad i=1,\dotsc, 6,
\label{eq:bcpo}
\end{equation} 
where $k_i \in \mathbb Z$ can be determined after an initial exploration of the solutions. Solutions to~\eqref{eq:rescaled} with the boundary conditions~\eqref{eq:bcpo} then correspond to periodic orbits in~$\T^6$. 

We need to supply initial data before continuing a periodic orbit in parameters. Namely, we need a periodic orbit that is a solution to the boundary value problem~\eqref{eq:rescaled}, \eqref{eq:bcpo}. The method described in \sref{sec:Param} is the key to providing reliable starting data for the continuation. This applies in the same way for computing an initial periodic orbit on~$\SDD$ and~$\SSD$. 
We first consider the reduced dynamics in normal form obtained from the parameterization method, so the corresponding rescaled system is 
\begin{equation}
v^{\prime}=Th(v,\lambda),
\label{eq:rescaled2}
\end{equation}  
where $h=f_{\delta}$ as in~\eqref{eq:f_delta} is the reduced dynamics in normal form~\eqref{eq:nf_SDD} or~\eqref{eq:nf_SSD} and $v= (v_1(t),v_2(t),v_3(t))\in \T^3$. 
The normal forms for~$\SDD$ and~$\SSD$ allow us to know where to look for initial conditions that converge to a stable periodic orbit on the corresponding invariant torus.
Once we set such an initial condition $v=v^0\in \T^3$, we solve the boundary value problem \eqref{eq:rescaled2} with
\begin{equation}
\begin{split}
v(0)&=v^0,\\
v(1)&=v^1,    
\end{split}
\label{eq:bc_endpoints}
\end{equation}
where $v^1\in \T^3$ is a vector of internal parameters. 
We also consider the boundary conditions
\begin{equation}
v_i(1)-v_i(0)-2l_i\pi=\eta_i, \quad i=1,2,3.
\label{eq:bc_po3}
\end{equation} 
Here $\eta=(\eta_1,\eta_2,\eta_3)$ is also a vector of internal parameters that is used to monitor when a solution to~\eqref{eq:rescaled2} with initial condition $v(0)=v^0$ is periodic in $\T^3$, for suitable fixed $(l_1,l_2,l_3)\in \mathbb Z^3$. 
We start the continuation of \eqref{eq:rescaled2}--\eqref{eq:bc_po3} setting $v^1=v^0$ and with the integration time~$T$,~$v^1$ and~$\eta$ free during the continuation, so that the problem is well-posed. 
After some initial exploration, we can detect values for $l_1,l_2$~and~$l_3$; for instance, on~$\SDD$ we have $l_1=2$, $l_2=l_3=-2$.
These values come from the analysis of the time series of the solution to the boundary value problem with the given initial condition, and are such that the solution behave like a periodic orbit. Since the orbit is stable in the reduced dynamics, this approximation is enough for AUTO to correct it and continue it in a the periodic orbit setting implemented in this paper. 
We stop the continuation when $\eta=0$, thus obtaining an approximation of a periodic orbit in~$\T^3$. 
Once we have this initial periodic orbit, we can embed it in~$\T^6$ as described in \sref{sec:dyn_tori} in order to get an initial periodic orbit that is a solution to \eqref{eq:rescaled}, \eqref{eq:bcpo} in~$\T^6$. 
The values of $k_1$--$k_6$ in~\eqref{eq:bcpo} are obtained through the embedding as well.

Once we have the starting data for a periodic orbit in~$\T^6$, we can continue it in system parameters. 
To look for bifurcations of these periodic orbits, we continue~\eqref{eq:rescaled}, \eqref{eq:bcpo} in~$\delta$. 
We add the integral phase condition 
\begin{equation}
\int_0^1\langle u(t), u'_{\tiny{\mbox{old}}}(t)\rangle dt=0,
\end{equation}
where~$u_{\mbox{\tiny{old}}}(t)$ is the previous solution computed in the continuation. 
This is to ensure uniqueness of the computed periodic orbit~\cite{DGK03}. 
We allow the coordinates of the end points~$u(0)$ and~$u(1)$ to move accordingly and can monitor them as internal parameters via extra user-defined boundary conditions. 
More importantly, in the {\sc Auto}-constants file we set the problem type as {\sc IPS}=7, so that we can monitor the Floquet multipliers of the periodic orbits obtained during the continuation in~$\delta$. 
One must be careful when using this option and monitor that one of the Floquet multipliers is always equal to one in order to avoid inaccuracy in the calculation.

Note that the full system on~$\T^6$ still has a continuous $\T$~symmetry that acts by a common phase shift to all oscillators. 
Since this leads to lack of hyperbolicity, we consider in practice the reduced system on~$\T^5$.
Specifically, we set $\theta_{\sigma,k}\to \theta_{\sigma,k}-\theta_{1,1}$, that is, we consider phase differences of each of the oscillators with respect to~$\theta_{1,1}$. 
This is the actual system numerically studied in \sref{sec:dyn_tori} with the setup explained above, with the corresponding changes considering its dimension. 
Here, in the boundary conditions \eqref{eq:bcpo} defining a periodic orbit, the information from~$k_1$ is passed to the rest of the constants via $k_i \to k_i - k_1$ for $i=2, \dotsc, 5$.

\end{document}